\documentclass[12pt, a4paper]{amsart}

\usepackage{amsmath, amsthm, amssymb}
\usepackage{graphicx}
\usepackage{enumerate}
\usepackage[top=2.2cm, bottom=2.1cm, left=2.1cm, right=2.2cm]{geometry}

\usepackage{stmaryrd}  

\newtheorem{definition}{Definition}
\newtheorem{lemma}{\bf Lemma}
\newtheorem{proposition}{\bf Proposition}
\newtheorem{theorem}{\bf Theorem}
\newtheorem{corollary}{\bf Corollary}

\renewenvironment{proof}{\noindent {\bf Proof: }}{\rm\\}
\theoremstyle{definition}
\newtheorem{remark}{Remark}{\rm}
\newtheorem{example}{Example}{\rm}

\usepackage{color}

\usepackage[T1]{fontenc}

\usepackage{algorithm, algpseudocode}

\makeatletter
\renewcommand{\p@algorithm}{\arabic{algorithm}\expandafter\@gobble}
\makeatother

\newcounter{step}[algorithm]
\newcommand\STEP[2][\(\triangleright\)]{%
	\refstepcounter{step}
	\vskip 0.25\baselineskip
	\item[]\hskip -\algorithmicindent #1 \textbf{Step \arabic{step}}%
	\ifthenelse{\equal{\unexpanded{#2}}{}}{}{ (\texttt{#2})}%
	\textbf{.}%
}
\newenvironment{algo}{\algo}{}
\def\algo#1\end{%
	\noindent\fbox{%
	\begin{minipage}[b]{\dimexpr\columnwidth-\algorithmicindent\relax}
	\begin{algorithmic}
	#1
	\end{algorithmic}
	\end{minipage}
	}%
\end}

\makeatletter
\renewcommand{\fnum@algorithm}{\fname@algorithm}
\makeatother

\usepackage{comment}
\usepackage{mathrsfs}

\usepackage{graphicx}
\DeclareFontEncoding{LS1}{}{}
\DeclareFontSubstitution{LS1}{stix}{m}{n}
\DeclareSymbolFont{stixletters}{LS1}{stix}{m}{it}
\DeclareMathAccent{\cev}{\mathord}{stixletters}{"91}
\DeclareMathAccent{\vec}{\mathord}{stixletters}{"92}
\DeclareMathAccent{\vecev}{\mathord}{stixletters}{"95}
\DeclareMathAccent{\cevstar}{\mathord}{stixletters}{"91}

\DeclareMathOperator*{\argmin}{argmin}

\makeatletter
\newsavebox\myboxA
\newsavebox\myboxB
\newlength\mylenA

\newcommand*\xoverline[2][0.75]{%
    \sbox{\myboxA}{$\m@th#2$}%
    \setbox\myboxB\null
    \ht\myboxB=\ht\myboxA%
    \dp\myboxB=\dp\myboxA%
    \wd\myboxB=#1\wd\myboxA
    \sbox\myboxB{$\m@th\overline{\copy\myboxB}$}
    \setlength\mylenA{\the\wd\myboxA}
    \addtolength\mylenA{-\the\wd\myboxB}%
    \ifdim\wd\myboxB<\wd\myboxA%
       \rlap{\hskip 0.5\mylenA\usebox\myboxB}{\usebox\myboxA}%
    \else
        \hskip -0.5\mylenA\rlap{\usebox\myboxA}{\hskip 0.5\mylenA\usebox\myboxB}%
    \fi}
\makeatother

\usepackage{tikz-cd}

\begin{document}

\title[Convergence of the EM algorithm]{Convergence of the EM algorithm via proximal techniques}
\author{Dominikus Noll}
\thanks{$^*$Institut de Math\'ematiques, Universit\'e de Toulouse, France}
\date{}

\begin{abstract}
We investigate convergence of the expectation maximization algorithm 
by representing it as a generalized proximal
method. Convergence of iterates and not just in value is investigated under natural hypotheses such as definability
of the incomplete data log-likelihood in the sense of o-minimal structure theory. 
\\

\\

\end{abstract}
\maketitle

\section{Introduction}
The EM algorithm assures monotone decrease of the incomplete data negative  log-likelihood \cite{dlr}, but
convergence of the iterates may fail in various ways \cite{wu}. 
Without coercivity iterates
may escape to infinity while values converge. Even when iterates stay bounded,
they may still fail to converge, cycle \cite{vaida}, or generate a continuum of accumulation points. 
Convergence analysis is further complicated
when iterates tend to 
the boundary of the natural parameter domain, where the likelihood is typically 
not well-behaved.

Even when  convergent, instead of approaching a global minimum of the incomplete data negative log-likelihood, EM
iterates may
 go to local minima, saddle points, or even to local maxima \cite{lachlan}. 
 This is a well-known phenomenon in non-convex optimization, where in practice it is usually
acceptable to find good local extrema.

In this work we allow  the maximum likelihood problem to include
parameter constraints. This yields a convenient way to model curved
exponential families,  but constraints
may also convey
prior knowledge about the unknown parameter, allow to implement restricted maximum likelihood \cite[p. 191]{lehmann},\cite{bentler,schoenberg97,phillips91,shi}, 
deal with truncated families \cite{levine},
keep iterates away from the boundary of the natural parameter domain, 
or simply  force boundedness, 
see e.g. \cite{liu_rubin1994,Hathaway83,Nettleton1999,Kim_Taylor1995}.  While  practical,
constraints further complicate convergence analysis of the EM algorithm.

There is a similarity between the EM algorithm and the proximal point method (PPM), 
which had been observed in
the contributions 
\cite{tseng,chretien1,chretien3,chretien4,chretien5}. 
The quadratic penalty term in PPM is replaced
by a regularizer 
based on the Kullback-Leibler information distance.  Since
convergence of PPM without convexity has been investigated \cite{artacho,attouch_bolte,pennanen,kaplan,spingarn,rock_pp_new}, 
some of these techniques, so the intention,  may carry over
to Kullback-Leibler regularizers.
Another compelling reason to investigate this link is the fact that
PPM can be seen as a special instance of EM when the latter  includes constraints. 

Presently we take a fresh look at this line, 
adding as a new element definability of the incomplete
data log-likelihood in the sense of o-minimal structure theory \cite{dries,dries_book,dries_miller,wilkie},
a hypothesis always met in practice. 
Our investigation reveals that the Kullback-Leibler distance has only a partial regularizing effect concentrated on a linear subspace, whose dimension depends on the rank of the 
conditional Fisher information matrix of missing data, given the observed datum.
In consequence, even under definability, only convergence of the projection of iterates on this subspace can be derived.
Further elements are needed to assure convergence of the full EM sequence. The positive aspect 
is
that this gives us clues as to why instances of EM fail to converge.

Additional insight into convergence of the full EM sequence $\theta^{(k)}$ is gained by
the point of view of \cite{bregman}, where the EM algorithm
had been interpreted as a method of
alternating Bregman projections between data and model sets. EM iterates $\theta^{(k)}$ generated in the M step arise in tandem with
iterates $\vartheta^{(k)}$ 
generated in the E step, and the algorithm alternates between these. Since the $\vartheta^{(k)}$ converge under fairly general hypotheses, also involving definability,
this adds new occasions to deduce convergence of the $\theta^{(k)}$.

While the EM algorithm is in general expected to converge linearly, cases of
sub-linear rates  have been reported, see  \cite[p.34]{lachlan}.
As part of our analysis we obtain a sublinear worst case rate $\|\theta^{(k)} - \theta^*\| = O(k^{-\rho})$
for some $0 < \rho < \infty$, which can be certified under quite natural assumptions.
The case where linear rates occur is also precisely delimited.

Recent approaches to convergence of the EM algorithm are \cite{kunstner,khan}, where the authors use mirror descent
to apply results from non-linear optimization, and  
\cite{caprio}, where  the Polyak-\L ojasiewicz inequality  allows the authors to derive
complexity results in a non-parametric setting. 
As the Polyak-\L ojasiewicz inequality is an instance of the
Kurdyka-\L ojasiewicz inequality, it is included in our present analysis, where
it gives the case of
linear convergence.

The structure of the paper is as follows. Section \ref{preparation} recalls facts from
optimization and definability theory. Section \ref{sect_EM} recalls the set-up
of the EM algorithm, including  the case of curved exponential families.
Section \ref{sect_KL} concerns Kullback-Leibler information, its role as a regularizer, and 
its link with Fisher information of missing data.
Parameter dimension reduction for
the conditional family follows in Section \ref{sect_dim}, revealing the partial character of the Kullback-Leibler regularizer.
Convergence under partial regularization is proved in Section \ref{sect_convergence},
a worst case rate given in Section \ref{sect_rate}.  Section \ref{convergence_EM} applies this to the EM algorithm,
followed by
cases where partial convergence can be upgraded to convergence of the full sequence $\theta^{(k)}$. 
Alternating Bregman projections come into play in Section \ref{sect_alter}. Extensions beyond the exponential family are discussed in Section \ref{sect_beyond}.
Examples are given in Section \ref{sect_examples}.

\section*{Notation}
For a function $\Psi(x,y)$ we write $\nabla_x \Psi$ of $\nabla_1\Psi$ for the derivative with respect to $x$, 
$\nabla_y\Psi$, $\nabla_2\Psi$ for the derivative with respect to $y$.
Second derivatives twice with respect to $x$ are $\nabla_{11}^2 \Psi(x,y)$ or $\nabla^2_{xx}\psi(x,y)$, and similarly
 $\nabla_{22}^2\Psi(x,y)$ or $\nabla^2_{yy} \Psi(x,y)$ and $\nabla^2_{xy}\Psi(x,y)$ for a mixed second derivative. The subdifferential is understood in the sense of \cite{mord,rock} and denoted $\partial \psi(x)$. For a function $\Psi(x,y)$
 we have $\partial_1 \Psi(x,y) = \partial \Psi(\cdot,y)(x)$, and $\partial_2\Psi(x,y)= \partial \Psi(x,\cdot)(y)$.
 The euclidean scalar product and norm on $\mathbb R^n$ are $x \cdot y$ and $\|x\|$. Euclidean balls are $B(x,\delta)$, and
 the euclidean $\delta$-neighborhood of a set $K$
 is $N(K,\delta)$.
 
\section{Preparation}
\label{preparation}
In this section we recall facts from optimization and definability theory. We follow the convention that iterates in general optimization algorithms
are denoted $x_k$, while when specifying to the EM algorithm iterates, being parameters to be estimated, will be termed $\theta^{(k)}$.

\subsection{Proximal point algorithm}
The classical proximal point method 
for a proper lower semi-continuous function $f:\mathbb R^n \to \mathbb R \cup \{\infty\}$
generates iterates $x_k$ via
\begin{equation}
\label{ppm}
x_{k+1} \in \argmin_{x \in \mathbb R^n} f(x) + \frac{1}{2\lambda_k} \|x - x_k\|^2,
\end{equation}
where $\lambda_k > 0$ and
where $\|\cdot\|$ is the euclidean norm \cite{martinet1,martinet2}. For convex $f$ it is known \cite{rock_pp,guler} that when $f$ has a minimum, 
the sequence $x_k$ converges to $x^*\in \argmin f$ 
iff $\sum_{j=1}^k \lambda_j \to \infty$ ($k\to \infty$),
and the speed is even super-linear
when $\lambda_k \to \infty$,  \cite{rock,guler}, \cite[Thm. 2.1]{luque}.
In the non-convex case
convergence is much harder to obtain, but partial results are known, cf. 
\cite{artacho,attouch_bolte,pennanen,kaplan,spingarn,rock_pp_new}.

\subsection{Bregman proximal method}
It has been proposed to replace
the quadratic penalty  $\frac{1}{2\lambda_k} \|x - x_k\|^2$ by a Bregman regularizer, 
$\lambda_k^{-1} D(x,x_k)$, 
which leads to the scheme
\begin{equation}
\label{PPB}
x_{k+1} \in \argmin_{x\in \mathbb R^n}
f(x) + \lambda_k^{-1}D(x,x_k),
\end{equation}
see \cite{teboulle,bdl}.
Here,
given a function $\psi$ of Legendre type \cite{BB_legendre},\cite{R_legendre}, the Bregman distance associated with $\psi$ is
\begin{equation}
\label{D}
D(x,y) = \left\{
\begin{array}{ll}
\psi(x) - \psi(y) - \nabla \psi(y)\cdot (x-y) & \mbox{ if $y \in {\rm int}({\rm dom} \,\psi)=:G$ } \\
+\infty & \mbox{ otherwise}
\end{array}
\right.
\end{equation}
with the proximal point method corresponding to the special case $\psi(x) = \frac{1}{2}\|x\|^2$. 

In the context of exponential families, the $\psi$ arise as cumulant generating functions, or log-normalizers. 
Legendreness of $\psi$  is called steepness \cite{brown}, and $\nabla^2 \psi \succ 0$ corresponds to minimality of the family. 
In the following
we assume throughout that $\psi$ is of class $C^2$.

\subsection{More general regularizers}
Expanding on (\ref{ppm}) and (\ref{PPB}), we next envisage schemes of the form
\begin{equation}
\label{general}
x_{k+1} \in \argmin_{x \in \mathbb R^n} f(x) + \lambda_k^{-1} \Psi(x,x_k), 
\end{equation}
with even more general regularizers $\Psi(x^+,x)$, allowing $D(x^+,x)$  as special cases.
Like $D(x^+,x)$, 
$\Psi(x^+,x)$ will only be partially defined.    %
Taking Bregman regularizers as paragon, we propose the following set-up:
\begin{enumerate}
    \item[(i)] There exists an open set $G$ with $G \times G \subset  {\rm dom} \nabla^2_{11}\Psi =
    {\rm dom} \nabla_1\Psi \subset {\rm dom}\,\Psi \subset \xoverline{G} \times G$, and
    $\Psi, \nabla_1 \Psi, \nabla^2_{11}\Psi$ are jointly continuous on their domains.
    \item[(ii)] $\Psi\geq 0$, and $\Psi(x,x)=0$ for $x\in G$.
\end{enumerate}
We call $\Psi$ {\it separating} if it satisfies the stronger property
\begin{enumerate}
\item[(ii')] $\Psi\geq 0$ and $\Psi(x^+,x)=0$ iff $x^+=x$.
\end{enumerate}
We say that $\Psi$ has a {\it lower norm bound} at $x\in G$ if there exist $\delta > 0$ and $m>0$ such that
\begin{enumerate}
\item[(iii)]
$\Psi(y,z) \geq m\|y-z\|^2$ for all $y,z\in B(x,\delta)$.
\end{enumerate}
We say that $\Psi$ has {\it pointwise lower norm bounds} on a set $K \subset G$ if (iii) holds for every $x\in K$.
While $\delta,m$ depend on $x$, a standard compactness argument (see Lemma \ref{lower}) shows that we can get the same $\delta,m$ 
for all $x\in K$ when $K\subset G$
is compact.  
Finally, in view of (i) we may without loss of generality assume that
dom$f \subset \xoverline{G}$.

\subsection{Partial regularizers}
\label{sect_partial}
Suppose ${\Psi}_V(v^+,v)$ satisfies (i), (ii), 
but only for elements $v^+,v$ of a linear subspace $V$. Then  $\Psi(x^+,x) = {\Psi}_V(Px^+,Px)$, with $P$
the orthogonal projection onto $V$, gives a regularizer on $\mathbb R^n$. We call such $\Psi$  partial regularizers,
because their effect is limited to $V$-components $Px$ of iterates $x$, leaving $V^\perp$-components unaffected.
In $x$-space a partial regularizer satisfies (i) and (ii), while (ii') or (iii) could at best be satisfied in $V$.

\subsection{Kurdyka-\L ojasiewicz inequality}
\label{definability}
The following definition is crucial for our approach.
\begin{definition}
{\rm (Kurdyka-\L ojasiewicz inequality)}.
{\rm
A lower semi-continuous function $f:\mathbb R^n \to \mathbb R \cup \{+\infty\}$ has the K\L-property at $\bar{x}\in  {\rm dom}(\partial f)$ if there exist
$\gamma,\eta > 0$, a neighborhood $U$ of $\bar{x}$,  and a continuous concave function $\phi:[0,\eta) \to \mathbb R_+$, called  {\it de-singularizing function},  such that
\begin{itemize}
\item[i.] $\phi(0)=0$,
\item[ii.] $\phi$ is of class $C^1$ on $(0,\eta)$,
\item[iii.] $\phi'(s) > 0$ for $s\in (0,\eta)$,
\item[iv.] For all $x \in U \cap \{x: f(\bar{x}) < f(x) < f(\bar{x})+\eta\}$ the K\L-inequality
\begin{equation}
\label{KL}
\phi'(f(x)-f(\bar{x})) {\rm dist}(0,\partial f(x)) \geq \gamma
\end{equation}
is satisfied, where $\partial f$ is the subdifferential of \cite{rock}. 
\end{itemize}
}
\end{definition}

\begin{remark}
We say that $f$ satisfies the 
\L ojasiewicz inequality when the de-singularizing function is $\phi'(s) = s^{-\theta}$ for some $\theta\in [\frac{1}{2},1)$, which means $\phi(s) = \frac{s^{1-\theta}}{1-\theta}$. 
The case $\theta=\frac{1}{2}$ is sometimes singled out under the name  Polyak-\L ojasiewicz inequality. 
\end{remark}

\begin{remark}
When $K$ is a compact set on which $f$ has constant value $f(\bar{x}) = f^*$ for all $\bar{x}\in K$, then (\ref{KL}) holds uniformly on a neighborhood $U$
of $K$. 
\end{remark}

\begin{remark}
For convergence via the K\L-property  see e.g.
\cite{absil,attouch_bolte,a_b_r_s,a_b_s,bolte2,noll_jota,gerchberg}. 
It is well-known that definability
in an o-minimal structure \cite{dries_book,dries,dries_miller}, for short {\it definability}, implies the K\L-inequality. See \cite{kurdyka}, and for non-smooth $f$,
\cite[Thm. 11]{b_d_l_s}, where it is shown that $\phi$ may be chosen concave. 
We use the o-minimal structure $\mathbb R_{an}$ of globally sub-analytic sets, but also the larger
$\mathbb R_{an,\exp}$, allowing exponential and logarithm,
see   \cite{dries_book,bierstone,coste,shiota,wilkie,miller}.
\end{remark}

\section{EM algorithm}
\label{sect_EM}
We consider a family of probability measures $\mathbb P_\theta$ with densities $p(x,\theta)$ with regard to a $\sigma$-finite base measure $\mu$ on the complete data space $X$, 
$d\mathbb P_\theta(x) = p(x,\theta) d\mu(x)$, where $\theta \in \Theta$ is the unknown  parameter 
we want to estimate by maximum likelihood. 
However, it is not $x\in X$ which is observed, but a random variable $y = h(x)$ in the incomplete data space $Y$, where $h:X \to Y$ is measurable and typically non-invertible. The density of observed data $y$ with regard to the marginal $\nu=\mu\circ h^{-1}$
is therefore
\begin{equation}
    \label{q}
q(y,\theta) = \int_{h^{-1}(y)} p(x,\theta) d\mu_y(x),
\end{equation}
where the family  $(\mu_y)_{y\in Y}$ is a disintegration of the measure $\mu$ with regard to the
marginal $\nu = \mu\circ h^{-1}$ on  $Y$, each $\mu_y$ concentrated on $h^{-1}(y) \subset X$.  
Here disintegration means
\begin{equation}
    \label{int_thm}
\int_X f(x) d\mu(x) = \int_Y \left[ \int_{h^{-1}(y)} f(x) d\mu_y(x) \right] d\nu(y)
\end{equation}
for $\mu$-integrable $f$.
Substituting $f = \chi_{h^{-1}(B)} p(\cdot,\theta)$ leads to the relation
$$
(\mathbb P_\theta\circ h^{-1})(B)=
\mathbb P_\theta(h^{-1}(B)) = \int_B \left[ \int_{h^{-1}(y)} p(x,\theta) d\mu_y(x) \right] d\nu(y)
$$
justifying (\ref{q}). This allows us to define, for every $y\in Y$, the conditional density
\begin{equation}
    \label{k}
k(x|y,\theta) = \frac{p(x,\theta)}{q(y,\theta)}, \; x\in h^{-1}(y),
\end{equation}
with regard to the measure $\mu_y$. For
hypotheses needed to establish the existence of a disintegration
 see for instance \cite{dys}.

Given an available sample $y$, maximum likelihood in incomplete data space $Y$
is the optimization program
\begin{equation}
\label{ML}
\hat{\theta}\in\argmin_{\theta \in M} -\log q(y,\theta),
\end{equation}
where  $M \subset {\Theta}$ is a set of model parameters admitted for optimization. 
At this stage the rationale of the EM algorithm
assumes that minimization (\ref{ML}) is cumbersome, and that it would be  preferable algorithmically  to perform maximum likelihood estimation
in complete data space
\begin{equation}
    \label{MLc}
\tilde{\theta} \in \argmin_{\theta \in M} -\log p(x,\theta).
\end{equation}
Since no sample $x$ is available, (\ref{MLc}) cannot be performed directly, and    
recourse is taken to the following
iterative procedure. Given a current guess $\theta^{(k)}\in M$ of the unknown parameter, one computes
{\it for every} $\theta\in M$, the conditional expectation of  $\log p(x,\theta)$, 
given $y$ and $\theta^{(k)}$:
$$
Q(\theta,\theta^{(k)}) := \mathbb E_{\theta^{(k)}}(\log p(x,\theta)|y),
$$
which (for fixed $y,\theta^{(k)}$) gives a function of $\theta$. This is called the E step, 
based on the formula
$$
\mathbb E_\theta(\phi(x,\theta')|y) = \int_{h^{-1}(y)} \phi(x,\theta') k(x|y,\theta) d\mu_y(x),
$$
with
$
k(x|y,\theta)$
given by (\ref{k}).
Once this is obtained, one performs
the M step
$$
\theta^{(k+1)} \in \argmin_{\theta \in M} -Q(\theta,\theta^{(k)}),
$$
which gives the new model parameter estimate $\theta^{(k+1)}\in M$. 
As is well-known, 
the EM algorithm reduces the negative log-likelihood $-\log q(y,\theta^{(k)})$
at each step, and this is not altered by optimizing over $\theta\in M$,
nor when only a local minimum is computed.
However, as mentioned before, monotone decrease in function value
does {\it not} assure convergence of the iterates $\theta^{(k)}$, and it is
convergence of the iterates we presently 
scrutinize.

\begin{algorithm}[h!]
\caption{EM algorithm}\label{algo_em_1}
\noindent\fbox{%
\begin{minipage}[b]{\dimexpr\textwidth-\algorithmicindent\relax}
\begin{algorithmic}
\STEP{E step} Given current model parameter estimate $\theta^{(k)}\in M$, make the function
$\theta \mapsto Q(\theta,\theta^{(k)}) = \mathbb E_{\theta^{(k)}}(\log p(x,\theta)|y)$ on $M$ available for optimization. 
\STEP{M step} Compute $\theta^{(k+1)} \in \argmin_{\theta \in M} -Q(\theta,\theta^{(k)})$.
Back to step 1.
	\end{algorithmic}
\end{minipage}
}%
\end{algorithm}

\begin{remark}
A special case of constraints  are  {\it curved families}
$M = \{\theta \in \Theta: \theta = \theta(u), u\in U\}$, with $\theta(u)$ a re-parametrization of $\theta$, but our approach
allows more general sets.
A  typical instance of $M$ is given in
Example \ref{ex_amari}.
\end{remark}

\subsection{Properties of exponential families}
The densities $p(x,\theta)$  on  $X$
form a $n$-dimensional exponential family with regard to the base measure $\mu$  if they are of the form
\begin{equation}
\label{exp_family}
p(x,\theta) = e^{\theta\cdot T(x) - \psi(\theta)},
\end{equation}
where $T(x)$ is the sufficient statistic,   $\theta \in \Theta = \{\theta\in \mathbb R^n: p(\cdot,\theta) \in L^1(X,d\mu)\}$ the natural
parameter, and $\psi(\theta)$ the log-normalizer
defined on $\Theta$, satisfying
\begin{equation}
\label{psi}
\psi(\theta) = \log \int_X e^{\theta\cdot T(x)} d\mu(x).
\end{equation}
The natural parameter space can also be written as $\Theta=\{\theta\in \mathbb R^n: \int_X e^{\theta\cdot T(x)} d\mu(x) < \infty\}$ and is a convex subset
of $\mathbb R^n$. There is no loss in generality in assuming that $\Theta$ is of full dimension $n$, as otherwise  a parameter reduction
leading to an equivalent representation (\ref{exp_family}) with lower dimension $m < n$ can be performed:
\begin{lemma}
\label{lemma1}
Suppose the natural parameter space $\Theta$ is contained in an affine subspace of dimension $m<n$. Then there exists
an equivalent representation of  $\mathbb P_\theta$ as a $m$-dimensional exponential family
$d\mathbb P_\theta(x) = {p}'(x,{\theta}')d{\mu}'(x) = e^{{\theta}' \cdot {T}'(x) - {\psi}'({\theta}')} d{\mu}'(x)$,
${\theta}' \in {\Theta}' \subset \mathbb R^m$, where now ${\rm dim}({\Theta}')=m$. If in {\rm (\ref{exp_family})} the statistic $T(x)$ is affinely independent
on $X$, then so is 
${T}'(x)$. 
\end{lemma}

\begin{proof}
Without loss of generality write the parameter as $\theta = (\theta_1,\theta_2)$ with $\theta_2 = A\theta_1 + a$ for a matrix $A$
of size $(n-m)\times m$  of rank $n-m$. Then
$d\mathbb P_\theta(x)=p(x,\theta) d\mu(x) = e^{ \theta_1\cdot T_1(x) + \theta_2 \cdot T_2(x) - \psi(\theta)} d\mu(x) = 
e^{ \theta_1 \cdot (T_1(x) + A^TT_2(x)) - {\psi}(\theta_1,A\theta_1+a)} e^{a\cdot T_2(x)} d\mu(x)
= e^{ \theta_1\cdot {T}'(x) - {\psi}'(\theta_1)} d{\mu}'(x)$,
with ${T}'(x) = T_1(x) + A^TT_2(x)$, ${\psi}'(\theta_1) = \psi(\theta_1,A\theta_1+a)$, and $d{\mu}'(x) = e^{a\cdot T_2(x)} d\mu(x)$, and where
the parameter space ${\Theta}' =\{\theta_1: (\theta_1,A\theta_1+a) \in \Theta\}$ is now of full dimension $m$. Since $m$ is smallest possible, there could no longer be
any affine dependence among the $\theta_1\in {\Theta}'$.
Note also that $\mu \ll {\mu}'$ and ${\mu}' \ll \mu$ gives equivalence of the two representations.

To conclude, suppose $a\cdot T(x)$ constant a.e. implies $a=0$. Then if ${a}' \cdot {T}'(x) = c$ for almost all $x\in X$, we have $({a}' , A^T{a}') \cdot (T_1(x),T_2(x)) = c$, hence
$({a}', A^T{a}')=(0,0)$, which gives ${a}'=0$, so that ${T}'$ is also affinely independent.
\hfill $\square$
\end{proof}

Assuming therefore that $\Theta$ has already full dimension $n$ in (\ref{exp_family}),
we denote the interior of $\Theta$ by $G$.
The function $\psi$ is also known as the cumulant generating function, because it satisfies
$$
\mathbb E_\theta[ T(x)] = \nabla \psi(\theta), \quad \mathbb V_\theta [T(x)] = \nabla^2 \psi(\theta),
$$
and similar relations for higher moments; cf. \cite{brown}. 

\begin{definition}
\label{def2}
{\rm 
A family (\ref{exp_family})  is called {\it minimal} if the functions $T_i(\cdot)$ are affinely independent, i.e., if
there exists no $a\not=0$
such that $a \cdot T(x) = const$ for $\mu$-a.a. $x\in X$. 
}
\end{definition}

\begin{lemma}
\label{lemma2}
Under minimality we have
$\nabla^2 \psi(\theta)  \succ 0$ for every $\theta\in G = {\rm int}(\Theta)$. Moreover, the mapping 
$\theta \mapsto p(\cdot,\theta)$ is injective on $G$.
\end{lemma}

\begin{proof}
1) From $\nabla^2 \psi(\theta) a=0$  follows
$0 = a\cdot \nabla^2 \psi(\theta) a = a\cdot \mathbb V_\theta[T(x)] a = \mathbb V_\theta [a\cdot T(x)] = \mathbb E_\theta|a\cdot T(x)- \mathbb E_\theta(a\cdot T(x))|^2$, hence
$a\cdot T(x) = \mathbb E_\theta[a\cdot T(x)] = a\cdot \nabla \psi(\theta)=const$ $\mu$-a.e., forcing $a=0$.  

2) Let
$p(\cdot,\theta)=p(\cdot,\theta')$ $\mu$-a.e., then
$\theta\cdot T(x)- \psi(\theta) =  \theta'\cdot T(x) - \psi(\theta')$ a.e., hence
$(\theta-\theta')\cdot T(x) + \psi(\theta')-\psi(\theta)=0$ a.e., so that
the vector
$a=\theta-\theta'$ renders $a\cdot T(x) = \psi(\theta)-\psi(\theta')$ constant a.e., forcing $\theta=\theta'$. 
\hfill $\square$
\end{proof}

An exponential family is called {\it steep} if the log-normalizer $\psi$ is of Legendre type \cite{rock,BB_legendre,brown}. The family is called {\it regular} if $\Theta$ is an open set, i.e., $\Theta = G$. In that case
the family is automatically steep, but the steep class is larger \cite{brown}.   

\subsection{EM algorithm for the exponential family}
\label{EM_exp}
The EM algorithm is
 sometimes characterized as going back and forth between  {\it completing the data} in the E step,  and
{\it maximum likelihood in complete data space} in the M step. This is not  true in general, cf. \cite{flury}, 
but holds when $p(x,\theta)$ belongs to an exponential family (\ref{exp_family}).
Namely, in that case,
$\log p(x,\theta) = \theta\cdot T(x) - \psi(\theta)$, hence  
$$
\mathbb E_{\theta^{(k)}}[\log p(x,\theta)| y] = 
\theta\cdot \mathbb E_{\theta^{(k)}}[T(x)| y]  - \psi(\theta),
$$
and the first term on the right {\it selects} a complete data statistic $t^{k+1}=T(x_{k+1})$, a fact
which one expresses by saying that
the E step consists in {\it completing the data}.  The
M step is unchanged, but due to the substitution of $T(x_{k+1})$,
leads to 
$Q(\theta,\theta^{(k)})=\log p(x_{k+1},\theta) =  \theta\cdot  T(x_{k+1}) - \psi(\theta)$, and may therefore
rightfully be referred to as
{\it maximum likelihood in complete data space}.
Altogether, for exponential families the EM algorithm has the form
\begin{algorithm}[h!]
\label{algo2}
\caption{EM algorithm for the exponential family}
\label{algo_em_2}
\noindent\fbox{%
\begin{minipage}[b]{\dimexpr\textwidth-\algorithmicindent\relax}
\begin{algorithmic}
\STEP{E step} Given current model parameter estimate $\theta^{(k)}\in M$, complete the data
by computing $T(x_{k+1})=\mathbb E_{\theta^{(k)}}[T(x)|y]$. 
\STEP{M step} Compute $\theta^{(k+1)} \in \argmin_{\theta \in M} \psi(\theta) -  \theta\cdot T(x_{k+1})$.
Back to step 1.
	\end{algorithmic}
\end{minipage}
}%
\end{algorithm}

\begin{remark}
    When the sufficient statistic  is affine, $T(x) = Ax + b$, the E step may even be based on computing
    $x_{k+1}=\mathbb E_{\theta^{(k)}}(x|y)$, as then $\mathbb E_{\theta^{(k)}}[T(x)|y]=\mathbb E_{\theta^{(k)}}(Ax+b|y)= A \mathbb E_{\theta^{(k)}}(x|y) + b$, 
    making the expression {\it completing the data} is even more suggestive.
\end{remark}

We conclude this section by remarking that when complete data are from an  exponential family $p(x,\theta)$ on $X$ as in (\ref{exp_family}), then the conditional densities
$k(x|y,\theta)$ also constitute, for given $y\in Y$, an exponential family on the sample space $h^{-1}(y)$ with regard to the base measure $\mu_y$ 
arising from the disintegration of $\mu$. This can be seen from
\begin{equation}
\label{cond_exp}
k(x|y,\theta) = \frac{p(x,\theta)}{q(y,\theta)} = \displaystyle\frac{e^{\theta\cdot T(x) -\psi(\theta)}}
{\int_{h^{-1}(y)}  e^{\theta\cdot T(x') -\psi(\theta)} d\mu_y(x')}
=: e^{\theta\cdot  T(x) -{\psi}_y(\theta)},
\end{equation}
obtained on putting
\begin{equation}
    \label{tilde_psi}
{\psi}_y(\theta) = \log \int_{h^{-1}(y)} e^{\theta\cdot T(x)}  d\mu_y(x),
\end{equation}
which parallels (\ref{psi}), and which we can write  as $d\mathbb P_\theta^{x|y} = k(x|y,\theta) d\mu_y(x)$.
For an exponential family  the measures $k(\cdot|y,\theta) d\mu_y$ are mutually equivalent,
i.e., $k(\cdot|y,\theta) d\mu_y \ll k(\cdot|y,\theta') d\mu_y$ for all $\theta,\theta'\in  \Theta$.
In a general setting, this may for practical considerations be added as a hypothesis (cf. \cite[Sect III B]{chretien1}).

\begin{remark}
Even when complete data $p(x,\theta)$ are from an exponential family,   
{\it this need not be true 
for incomplete data}  $q(y,\theta)$.
The exponential structure may be lost just because data are missing.
But there are also cases where specific families $q(y,\theta)$ {\it not} of exponential type are deliberately arranged as {\it missing data from an exponential family}. 
This terminology goes back to Sundberg \cite{sundberg},  who gives a variety of examples $q(y,\theta)$, including finite mixtures of exponential families, 
censored data, convolutions, folded distributions,
the negative binomial distribution, and much else.
Presently we extend this to include parameter constraints $\theta \in M$, so that complete data $p(x,\theta)$
from which $q(y,\theta)$ are derived may e.g.  be curved.
\end{remark}

Concerning the well-posedness of the EM sequence, we have to bear in mind that $M$ has to be a closed set,
because the objective $f$ has to be lower semi-continuous. This may be in conflict with the fact that $\Theta$
is in general not closed. We may therefore only assume $M \subset \xoverline{\Theta}$, so that $M \cap \partial \Theta$ may be non-empty,
and points in this set cause trouble.

\section{Kullback-Leibler information measure}
\label{sect_KL}
The Kullback-Leibler information distance on $X$ is defined as
$$\mathcal K(q||p) =\mathbb E_q\left(\log \frac{q}{p} \right)= \int_X q(x) \log \frac{q(x)}{p(x)} d\mu(x),$$
and in the parameter-dependent case we use the notation
$$
K(\theta||\theta^+) = \mathcal K(p(\cdot,\theta)||p(\cdot,\theta^+)).
$$
When restricted
to $h^{-1}(y)$ the KL-distance takes the form
\begin{equation}
    \label{KLy}
K_y(\theta||\theta^+) = \mathbb E_\theta\left(\log \frac{k(\cdot|y,\theta)}{k(\cdot|y,\theta^+)}\bigg| y\right)
=
\int_{h^{-1}(y)} k(x|y,\theta) \log \frac{k(x|y,\theta)}{k(x|y,\theta^+)} d\mu_y(x),
\end{equation}
where $\mu_y$ arises from the disintegration of $\mu$, and where the value is finite due to the hypothesis
$k(\cdot|y,\theta) d\mu_y \ll k(\cdot|y,\theta') d\mu_y$ for all $\theta,\theta'\in \Theta$.
The following is now a crucial observation.

\begin{proposition}
\label{prop1}
{\rm (See \cite[Prop. 1]{chretien1})}.
The {\rm EM} algorithm is a realization of the general scheme  {\rm (\ref{general})} 
with
$f(\theta) = -\log q(y,\theta) + i_{M}(\theta)$, $\lambda_k = 1$, and
$\Psi(\theta,\theta^{(k)}) = K_y(\theta^{(k)}||\theta	)$ given by {\rm (\ref{KLy})}.
\end{proposition}

\begin{proof}
Re-arranging (\ref{k}) and integrating, we have
$$
\log q(y,\theta) = \mathbb E_{\theta^{(k)}}(\log p(x,\theta)|y) - \mathbb E_{\theta^{(k)}}(\log k(x|y,\theta) | y)
$$
for arbitrary $k$. 
Hence
the M step in the EM algorithm is
$$
\theta^{(k+1)} \in \argmin_{\theta \in \mathbb R^n} -\log q(y,\theta) +i_M(\theta) - \mathbb E_{\theta^{(k)}} \left(\log k(x|y,\theta) \big| y\right).
$$
Adding the constant term $\mathbb E_{\theta^{(k)}}\left( \log k(x|y,\theta^{(k)}) \big| y\right)$ to the objective
does not change the optimization program, hence we have
\begin{equation}
\label{missing}
\theta^{(k+1)} \in \argmin_{\theta \in \mathbb R^n} -\log q(y,\theta) + i_{M}(\theta) - \mathbb E_{\theta^{(k)}}\left( \log \frac{k(x|y,\theta)}{k(x|y,\theta^{(k)})} \bigg| y  \right),
\end{equation}
and the last term equals $K_y(\theta^{(k)}||\theta)$.
\hfill $\square$
\end{proof}

\begin{remark}
1) This explains why the EM algorithm decreases the negative log-likelihood
$-\log q(y,\theta)$, as this is a general property of the scheme (\ref{general}), see also Theorem \ref{main}.

2)
    Kullback-Leibler divergence is separating in  function space, i.e.,
    $\mathcal K(q||p)=0$  implies $q=p$ a.e. Hence $K_y(\theta||\theta^+)=0$
    implies $k(\cdot|y,\theta) = k(\cdot|y,\theta^+)$ $\mu_y$-a.e. However, the latter does not always give $\theta = \theta^+$, 
    because the family $k(\cdot|y,\theta)$ is not necessarily minimal on $h^{-1}(y)$. 
\end{remark}

\subsection{Interpretation for the exponential family}
\label{sect_KL_exp}
We now specify the Kullback-Leibler divergence to the case of an exponential family.

\begin{proposition}
\label{Kull_Breg}
The Kullback-Leibler divergence of two distributions $p(x,\theta)$ and $p(x,\theta')$ belonging to the same exponential family is
$K(\theta||\theta^+) = D(\theta^+,\theta)$, where $D$ is the Bregman divergence induced by the log-normalizer $\psi$.
\end{proposition}

\begin{proof}
From (\ref{exp_family}), and since $\int_Xp(x,\theta) d\mu(x)=1$,  we have
$$
\psi(\theta) = \log \int_{X} e^{ \theta\cdot T(x)} d\mu(x).
$$
Differentiation with respect to $\theta$ gives
$$
\nabla \psi(\theta) = {\int_X T(x) e^{\theta\cdot T(x)} d\mu(x) }\bigg/{\int_{X} e^{\theta\cdot T(x)}d\mu(x)}.
$$
Now $e^{\psi(\theta)} = \int_X e^{\theta\cdot T(x)} d\mu(x)$, hence $\nabla \psi(\theta) =\int_X T(x) e^{\theta\cdot T(x) - \psi(\theta)  } d\mu(x)
= \int_X T(x) p(x,\theta) d\mu(x)=E_\theta [T(x)]$, the expectation  of the random variable $T(x)$ with respect to the distribution
$d\mathbb P_\theta=p(\cdot,\theta) d\mu$ (see also \cite[(1.57)]{lachlan}). Then
\begin{align}
\label{reversed}
\begin{split}
K(\theta||\theta^+) &= \int_X p(x,\theta) \log\frac{p(x,\theta)}{p(x,\theta^+)}\, d\mu(x) \\
&= \int_X p(x,\theta) \left[ \psi(\theta^+)-\psi(\theta) + (\theta-\theta^+)\cdot T(x) \right] d\mu(x)\\
&=\int_X p(x,\theta) \left[  D(\theta^+,\theta) + (\theta^+-\theta)\cdot \nabla \psi(\theta) + (\theta-\theta^+)\cdot T(x) \right] d\mu(x)\\
&=D(\theta^+,\theta) + \int_X p(x,\theta) \left[ (\theta^+-\theta) \cdot  (\nabla \psi(\theta)-T(x))\right] d\mu(x)\\
&= D(\theta^+,\theta) + (\theta^+-\theta) \cdot  (\nabla \psi(\theta) - E_{\theta}[T(x)]) \\
&= D(\theta^+,\theta).
\end{split}
\end{align}
This proves the claim.
\hfill $\square$
\end{proof}

Bregman distances or divergences are usually considered for functions $\psi$ of Legendre type \cite{BB_legendre,rock}. As already mentioned, 
for log-normalizers this is called
steepness \cite[Def. 3.2]{brown}. Most exponential families in practice are regular, i.e., $\Theta$ is open, in which case steepness follows automatically.

\begin{lemma}
Suppose the exponential family $p(x,\theta)$ is minimal. Then the Bregman divergence induced by the log-normalizer $\psi$ is separating,
i.e., $D(\theta^+,\theta)=0$ implies $\theta^+=\theta$.
\end{lemma}

\begin{proof}
From $D(\theta^+,\theta)=0$ we get
$\psi(\theta^+)-\psi(\theta) - \nabla \psi(\theta)\cdot (\theta^+-\theta)=0$. Taylor expansion at $\theta$ gives
$\psi(\theta^+) = \psi(\theta) + \nabla \psi(\theta)\cdot(\theta^+-\theta) + \frac{1}{2} (\theta^+-\theta)\cdot  \nabla^2 \psi(\bar{\theta}) (\theta^+-\theta)$ for some
$\bar{\theta}$ on the open segment joining $\theta^+$ and $\theta$, and depending on $\theta,\theta^+$. Hence
$(\theta^+-\theta)\cdot \nabla^2 \psi(\bar{\theta}) (\theta^+-\theta) = 0$. But minimality gives $\nabla^2 \psi(\bar{\theta}) \succ 0$ by Lemma \ref{lemma2}, hence $\theta^+=\theta$.
\hfill $\square$
\end{proof}

\begin{lemma}
Suppose the exponential family $p(x,\theta)$ is minimal. Then $\nabla \psi$ is injective on $G$.
For $\eta = \nabla \psi(\theta) \in G^* = {\rm int}({\rm dom} \,\psi^*)$ we have $\nabla\psi^*(\eta) = \theta$.
In particular, if the family is steep, then
$\nabla \psi^*$ is the inverse of $\nabla\psi$, with $G = {\rm int}({\rm dom}\, \psi)$ mapped {\rm 1-1} onto
$G^* = {\rm int}({\rm dom}\, \psi^*)$.
\end{lemma}

\begin{proof}
Let $\nabla \psi(\theta) = \nabla \psi(\theta^+)$. For a test vector $h$, Taylor expansion of $\theta \mapsto h \cdot \nabla \psi(\theta)$ at $\theta^+$
gives
$h\cdot \nabla \psi(\theta) = h\cdot \nabla \psi(\theta^+) + h\cdot \nabla^2 \psi(\bar{\theta}) (\theta-\theta^+)$ for some $\bar{\theta}=\bar{\theta}(\theta,\theta^+,h)$
on the open segment joining $\theta$ and $\theta^+$ and depending on $\theta,\theta^+,h$. Matching this with the first equation above shows
$h \cdot \nabla^2\psi(\bar{\theta})(\theta-\theta^+)=0$. Taking as test vector $h=\theta-\theta^+$ gives
$(\theta-\theta^+) \cdot \nabla^2 \psi(\bar{\theta})(\theta-\theta^+)=0$ for $\bar{\theta} = \bar{\theta}(\theta,\theta^+,\theta-\theta^+)$,
and since $\nabla^2 \psi(\bar{\theta}(\theta,\theta^+,\theta-\theta^+)) \succ 0$, we have $\theta = \theta^+$.

For the second part,  recall that $\partial \psi$, $\partial \psi^*$ are inverses of each other
in the sense
$\eta \in \partial \psi(\theta)$ iff $\theta\in \partial \psi^*(\eta)$; cf.  \cite[Cor. 23.5.1]{rock}. 
By
strict convexity of $\psi$ its conjugate $\psi^*$ is differentiable on $G^*={\rm int \, dom} \,\psi^*$. 
Hence if $\eta = \nabla \psi(\theta)\in G^*$, then $\nabla \psi^*(\eta) = \theta$.
Since we may have $\nabla \psi(G) \not\subset G^*$,  all we know
about  $\eta = \nabla\psi(\theta) \in \partial G^*$ is
$\theta \in \partial \psi^*(\eta)$.
When $p(x,\theta)$ is steep, then $\partial \psi(\theta) = \emptyset$ for $\theta\in \partial \Theta$, and then
$\nabla \psi$ maps $G$ 1-1   into $G^*$ with inverse $(\nabla \psi)^{-1}=\nabla \psi^*$. 
\hfill $\square$
\end{proof}

\begin{remark}
In general one has $G={\rm int} ({\rm dom} \,\psi)\subset {\rm dom} (\nabla \psi) \subset {\rm dom} \,\psi = \Theta$,
and $G^* ={\rm int} ({\rm dom} \,\psi^*)\subset {\rm dom} (\nabla \psi^*) \subset {\rm dom} \,\psi^*$;
cf. \cite[Thm. 23.4]{R_legendre}.
 \end{remark}

\subsection{Fisher information of missing data}
\label{sect_fisher}
The regularizer in (\ref{missing}) has a statistical interpretation.
It is well-known that
$K_y(\theta||\theta)=0$, and also 
$K_y \geq 0$,  hence for fixed $\theta$,  the global minimum $0$ of 
$\theta^+ \mapsto K_y(\theta||\theta^+)$ is attained in particular at $\theta^+=\theta$. Then clearly $\nabla_2 K_y(\theta||\theta)=0$ and
$\nabla^2_{22}K_y(\theta||\theta) \succeq 0$ from the necessary optimality conditions. We investigate whether we may
expect the stronger sufficient optimality condition
$\nabla^2_{22}K_y(\theta||\theta) \succ 0$. 
Going back to the definition, we have 
$$
K_y(\theta||\theta^+) = \int_{h^{-1}(y)} k(x|y,\theta) \log \frac{k(x|y,\theta)}{k(x|y,\theta^+)} d\mu_y(x)  =\mathbb E_\theta\left[ \log \frac{k(\cdot|y,\theta)}{k(\cdot|y,\theta^+)} \bigg| y\right].
$$
Differentiating twice with respect to $\theta^+$ (cf. \cite[Thm. 5.8, sect. 7.1]{barndorff}) gives
$$
\nabla^2_{22} K_y(\theta||\theta^+) = \mathbb E_\theta[-\nabla^2_{\theta^+\theta^+} \log k(\cdot|y,\theta^+) | y],
$$
hence we obtain
$$
\nabla^2_{22} K_y(\theta||\theta^+)\big|_{\theta^+=\theta} = \mathbb E_\theta[-\nabla^2_{\theta\theta} \log k(\cdot|y,\theta) | y] =: \mathcal I_m(\theta,y),
$$
which  is recognized as the conditional expected Fisher information matrix of missing data, given $y$; cf. \cite[3.52]{lachlan}.  According to the missing information principle
\cite{orchard}, $\mathcal I_m(\theta,y)$ gives the expected loss of information between complete and incomplete data.
Differentiating the identity $\mathbb E_\theta[ \nabla_\theta\log k(\cdot|y,\theta) | y] = 0$   with respect to $\theta$ (see again  \cite[Thm. 5.8, sect. 7.1]{barndorff})
gives the alternative formula
\begin{equation}
\label{stimmts}
\mathcal I_m(\theta,y) = 
\mathbb E_\theta[ \nabla_\theta k(\cdot|y,\theta) \nabla_\theta k(\cdot|y,\theta)^T | y] 
\succeq 0.      
\end{equation}
A bit more can be said in the case of an exponential family.

\begin{lemma}
\label{lemma5}
    Suppose $k(x|y,\theta)$ is an exponential family. Then
    $\nabla^2_{22} K_y(\theta||\theta) = \nabla^2 {\psi}_y(\theta)$. In addition, if the family is minimal with regard to the sample space $h^{-1}(y)$, then $\nabla^2 \psi_y(\theta)\succ 0$.
\end{lemma}

\begin{proof}
       From (\ref{reversed}) we  get $K_y(\theta||\theta^+) = D_y(\theta^+,\theta)$ for the Bregman distance $D_y$ induced by $\psi_y$,
    and then $\nabla^2_{\theta^+\theta^+} K_y(\theta||\theta^+) = \nabla^2 \psi_y(\theta^+)$. The last part follows with Lemma \ref{lemma2}.
    \hfill $\square$
     \end{proof}

\section{Dimension reduction for the conditional family}
\label{sect_dim}
We consider constrained maximum likelihood
with incomplete data from a $n$-dimensional exponential family, i.e., ${\rm im}(T) \subset \mathbb R^n$, where
dim$(\Theta)=n$. In view of Lemma \ref{lemma1},  we also assume that the complete data family $p(x,\theta)$ is minimal. However, this does
not mean that the conditional $n$-dimensional exponential family
$k(x|y,\theta)$ is also minimal on $h^{-1}(y)$. 
In fact,  the missing data case (Example \ref{ex_missing})
shows that we should rather expect the opposite.  
This calls for a dimension reduction argument.

\begin{proposition}
\label{prop4}
The $n$-dimensional conditional exponential family $d\mathbb P_\theta^{x|y} = k(\cdot|y,\theta) d\mu_y$  has 
an equivalent minimal representation as  $m$-dimensional exponential family
$d\mathbb P_{{\theta'}}^{x|y} = \bar{k}(x|y,{\theta}') d{\mu}_y$, where ${\theta}'=P\theta$ for an orthogonal projection $P$ on a $m$-dimensional subspace
$V$ of $\mathbb R^n$. The conditional distributions satisfy $\mathbb P_\theta^{x|y} = \mathbb P_{{\theta}'}^{x|y}$, with sufficient statistic ${T'} = P \circ T$
now minimal. 
Under this reduction the Kullback-Leibler divergences are related as
$K_y(\theta||\theta^+)= \bar{K}_y({\theta'}||{\theta'}^+)$. 
\end{proposition}

\begin{proof}
We have $T:X \to  \mathbb R^n$, and
$T(h^{-1}(y)) \subset t_0 + V$ for a linear subspace $V$ of $\mathbb R^n$ of minimal dimension $m$, where we may assume without loss of generality that $t_0\in V^\perp$.
Let $P$ be the orthogonal projection onto $V$.
Write $T(x) = t_0 + v(x)$ for $v(x)\in V$.
We have
\begin{align}
\label{19}
\begin{split}
\psi_y(\theta) &= \log \int_{h^{-1}(y)} e^{\theta \cdot t_0} e^{\theta \cdot v(x)} d\mu_y(x)\\
&= \theta \cdot t_0 + \log \int_{h^{-1}(y)} e^{P\theta \cdot v(x)} d\mu_y(x) 
= :\theta \cdot t_0 + \bar{\psi}_y(P\theta),
\end{split}
\end{align}
where we use the fact that $\theta \cdot v = P\theta \cdot v$ for $v\in V$, and where in consequence the rightmost term $\bar{\psi}_y(P\theta)$ depends only on $P\theta$.
Now $Pt_0=0$ implies $PT(x) = v(x)$, hence
\begin{align}
\label{20}
\begin{split}
\log k(x|y,\theta) &= \theta \cdot T(x) - \psi_y(\theta) \\
&= \theta \cdot t_0 + \theta \cdot v(x) - \psi_y(\theta)\\
&= \theta \cdot t_0 + \theta \cdot v(x) - \bar{\psi}_y(P\theta) - \theta \cdot t_0    \hspace{1cm}\mbox{ (using (\ref{19}))}\\
&= P\theta \cdot v(x) - \bar{\psi}_y(P\theta) =: \log \bar{k}(x|y,P\theta),
\end{split}
\end{align}
using again $\theta	 \cdot v = P\theta \cdot v$ for $v\in V$.
This gives the representation
\begin{equation}
\label{21}
d\mathbb P_\theta^{x|y} = k(\cdot|y,\theta) d\mu_y = \bar{k}(\cdot|y,{\theta}') d\mu_y = d\mathbb P_{{\theta}'}^{x|y},
\end{equation}
where the new sufficient statistic is $v =P\circ T$, and the new parameter is  ${\theta}' = P\theta$. 
Note that $v=P\circ T$ is affinely
independent on $h^{-1}(y)$, because
$T(h^{-1}(y)) - t_0 \subset V$ has full dimension $m$ in $V$ by the choice of $V$, and
we have $T(h^{-1}(y)) - t_0 = P[ T(h^{-1}(y)) - t_0] = v(h^{-1}(y))$, so that the latter has also full dimension in $V$.

It remains to compare the Kullback-Leibler divergences generated by both representations. 
\begin{align*}
\bar{K}_y({\theta'}|| {\theta'}^+) &= \int_{h^{-1}(y)} \log \frac{\bar{k}(x|y,{\theta'})}{\bar{k}(x|y,{\theta'}^+)} \bar{k}(x|y,{\theta'})  d\mu_y(x) \\
&=  \int_{h^{-1}(y)}  \log \frac{\bar{k}(x|y,{\theta'})}{\bar{k}(x|y,{\theta'}^+)} k(x|y,\theta)d\mu_y(x)    \hspace{1cm} \mbox{ (using (\ref{21}))}\\
&= \int_{h^{-1}(y)} \left[  (P\theta-P\theta^+)\cdot v(x) - \bar{\psi}_y(P\theta)+\bar{\psi}_y(P\theta^+) \right]  k(x|y,\theta)d\mu_y(x) \\
&= \int_{h^{-1}(y)}   \left[  (\theta-\theta^+)\cdot v(x) - {\psi}_y(\theta)+ \theta\cdot t_0+{\psi}_y(\theta^+) -\theta^+\cdot t_0\right]  k(x|y,\theta)d\mu_y(x) \\
&=  \int_{h^{-1}(y)}   \left[  (\theta-\theta^+)\cdot (v(x)+t_0) - {\psi}_y(\theta) +{\psi}_y(\theta^+) \right]  k(x|y,\theta)d\mu_y(x) \\
&=  \int_{h^{-1}(y)}   \left[  (\theta-\theta^+)\cdot T(x) - {\psi}_y(\theta) +{\psi}_y(\theta^+) \right]  k(x|y,\theta)d\mu_y(x) \\
&= \int_{h^{-1}(y)} \log \frac{k(x|y,\theta)}{k(x|y,\theta^+)}  k(x|y,\theta)d\mu_y(x) = K_y(\theta||\theta^+).\end{align*}
Finally, for notational beauty we write $v$ as $T'$.
\hfill $\square$
\end{proof}

The relationship between $\nabla^2_{\theta\theta} \psi_y(\theta)$ and $\nabla^2_{{\theta'}{\theta'}}\bar{\psi}_y({\theta'})$ is as follows.

\begin{corollary}
There exists a
$n\times n$ orthogonal matrix $Q$ with
\begin{equation}
\label{Q}
Q \nabla^2_{\theta\theta}\psi_y(\theta) Q^T = \left[\begin{array}{c|c} \nabla^2_{{\theta'}{\theta'}}\bar{\psi}_y({\theta'}) & 0 \\ \hline
0 & 0 \end{array}\right], \quad Q \theta = \left[ \begin{array}{c} {\theta'}\\\ \!\!{{\theta''}}\end{array}\right], \quad  \nabla^2_{{\theta'}{\theta'}}\bar{\psi}_y({\theta'}) \succ 0,
\end{equation}
and in particular $m={\rm dim}({\theta'})$ is the rank of $\nabla^2_{\theta\theta}\psi_y(\theta)$.
The KL-divergence ${K}_y({\theta}||{\theta}^+)$ depends only on the coordinate ${\theta}'$
and is separating on the subspace $V={\rm im}(P)$, where $P$ is the orthogonal projection 
 ${\theta}'=P\theta$. Moreover $\nabla^2_{22}\bar{K}_y({\theta'}||{\theta'}) \succ 0$.
\end{corollary}

\begin{proof}
The projections $P:\mathbb R^n \to V$ and $I-P:\mathbb R^n \to V^\perp$ used in the proof of Proposition \ref{prop4} are represented by
an orthogonal $n\times n$ matrix $Q$ giving  the change of variables in (\ref{Q}), where
$\theta'\in V$, $\theta''\in V^\perp$. Then $\bar{k}$ and the log-normalizer $\bar{\psi}_y$
depend only on ${\theta'}$. 

Since the family $\bar{k}(x|y,{\theta'})$ is minimal, we have indeed $\nabla^2_{{\theta'}{\theta'}} \bar{\psi}_y({\theta'}) \succ 0$ by
Lemma \ref{lemma2}. 
Since in a minimal family ${\theta'} \mapsto \bar{k}(\cdot|y,{\theta'})$ is 1-1, 
$\bar{K}_y({\theta'}||{\theta'}^+)=0$ gives in the first place $\bar{k}(\cdot|y,{\theta'})=\bar{k}(\cdot|y,{\theta'}^+)$
${\mu}_y$-a.e., and then ${\theta'}={\theta'}^+$. 
The last claim follows from Lemma \ref{lemma5}.
\hfill $\square$
\end{proof}

\begin{remark}
The rank $m$  depends on $y$. Since dim$({\theta'})=m$, we call ${\theta'}$ the {\it accurate} parameter,
unique up to an orthogonal change of coordinates. Its meaning is that the family $k(x|y,\theta)$ is overparametrized by the $n-m$ {\it spare} parameters $\theta''$,
and has a statistic with $n-m$ too many components $T_j(x)$,
while $\theta'$ maintains only the accurate number $m$ of parameters needed. 
\end{remark}

\begin{remark}
While 
the topological dimension of $T(h^{-1}(y))$ is typically smaller than $n$,  it is  possible that $m=n$ for the affine dimension of $T(h^{-1}(y))$. For instance, $T(h^{-1}(y))$ might be a space curve
in 3d-space, which has topological dimension 1, but affine dimension 3. This has the ironic consequence that our convergence result for such curved fibers
$h^{-1}(y)$  is {\it a priori} better than for the more likely case
where fibers are affine subspaces.
\end{remark}

What we have found is that $\Psi(\theta^+,\theta) = K_y(\theta||\theta^+)$ in (\ref{missing}) is a partial regularizer
on the subspace $V = {\rm im}(P)$ in the sense of Section \ref{sect_partial}. This calls now for our central convergence
result under partial regularization, which we give in the next section.

\section{Convergence with interiority}
\label{sect_convergence}

We prove partial convergence of the generalized proximal method (\ref{general}) under the assumption that 
the sequence of iterates is bounded and together with its set of accumulation points stays in the interior $G$ of the domain
of $\Psi$. Since the results are of general nature, we switch to the notation familiar in optimization.

We consider an extension where (\ref{general}) is solved approximatively
in the sense that
\begin{equation}
\label{error}
e_k = g_k + \lambda_{k-1}^{-1} \nabla_1 \Psi(x_k,x_{k-1}) \mbox{ and } f(x_{k}) + \lambda_{k-1}^{-1} \Psi(x_{k},x_{k-1}) \leq f(x_{k-1}),
\end{equation}
with $g_k\in \partial f(x_k)$
and a subgradient error  $e_k$ satisfying one of the following conditions:
\begin{itemize}
\item[a.] $\sum_k \lambda_{k-1}\|e_k\| < \infty$;
\item[b.] $\lambda_{k-1}\|e_k\| \leq M'  \|\nabla_1 \Psi(x_k,x_{k-1})\|$ for some fixed big $M'>0$;
\item[c.] $\|e_k\| \leq M'' \|g_k\|$ for some big $M''>0$.
\end{itemize}

\begin{theorem} {\bf (Partial convergence)}.
\label{main}
Suppose $f$ satisfies the K\L-inequality on $G$. Consider a sequence $x_k$ generated by the approximate proximal method, which is bounded and together with its accumulation points stays in $G$.
Let $\Psi_V$ have pointwise lower norm bounds and be separating on a subspace $V$ with projection $P$.
Assume $\lambda_k/\lambda_{k-1} \leq r < \infty$ and $\lambda_k \leq R < \infty$. Then the sequence $Px_k$ converges. When $\sum_k \lambda_k = \infty$, then 
at least one accumulation point of the $x_k$ is critical, and if $\lambda_k \geq \eta > 0$, then all accumulation points are critical.
\end{theorem}

\begin{proof}
1)
From (\ref{general}) in the case where $e_{k+1}=0$, respectively from (\ref{error}), we have
\begin{equation}
\label{descent2}
f(x_{k+1})+\lambda_k^{-1} \Psi(x_{k+1},x_k) \leq f(x_k) + \lambda_k^{-1} \Psi(x_{k},x_k) = f(x_k),
\end{equation}
using $\Psi(x,x)=0$. For $\Psi(x_{k+1},x_k) > 0$ values  are strictly decreasing. The case $f(x_{k+1}) = f(x_k)$ only occurs when
$\Psi(x_{k+1},x_k)=0$, which by (\ref{error}) implies $e_{k+1}\in \partial f(x_{k+1})$. Here the algorithm stops when $e_{k+1}=0\in \partial f(x_{k+1})$,
but continues for $e_{k+1} \not=0$.
We may therefore
concentrate on the case where the algorithm does not terminate finitely.

2)
By hypothesis the sequences $x_k$ and $f(x_k)$ are bounded, accumulation points of the $x_k$ are in $G$,
and by monotone convergence we have
$f(x_k) \to f^*\in \mathbb R$.
From this and (\ref{descent2}) we immediately get $\lambda_k^{-1} \Psi(x_{k+1},x_k) \to 0$, using $\Psi \geq 0$. 
Since $\lambda_k^{-1} \geq R^{-1} > 0$, we deduce
$\Psi(x_{k+1},x_k) \to 0$. We argue that this implies $P(x_{k+1}-x_k) \to 0$. Indeed, assume that there is an infinite subsequence $k\in N\subset \mathbb N$
with $\|P(x_{k+1}-x_k)\|\geq \epsilon > 0$.
Using boundedness of $x_k$, extract a sub-subsequence $k'\in N'\subset N$ such that
$x_{k'}\to x$, $x_{k'+1}\to x'$. Then $\Psi(x_{k'+1},x_{k'}) \to \Psi(x',x)=0$, and since $\Psi$ is separating on $V$,
we have $Px'=Px$, forcing $P(x_{k'+1}-x_{k'}) \to 0$, a contradiction.

3)
We argue that this implies $\Psi(x_{k+1},x_k) \geq m\|P(x_{k+1}-x_k)\|^2$ from some counter onwards. Indeed, let $K$ be the 
compact set of accumulation points of the $x_k$. Applying Lemma \ref{lower} to the regularizer $\Psi_V$ in the space $V$,
we find that there exist $m,\delta > 0$ such that
$m\|P(y-z)\|^2 \leq \Psi(y,z)$ for all $y,z\in N(K,\delta)$ with $\|P(y-z)\| < \delta$. But clearly $x_{k+1},x_k\in N(K,\delta)$ from
some counter $k$ onward, and as $P(x_{k+1}-x_k) \to 0$, we also have $\|P(x_{k+1}-x_k)\| < \delta$ from some counter onwards.
This proves the claim.

4)
The sequence $x_k$ being bounded, its set of accumulation points $K$ is compact, and $f$ has constant value $f^*$
on $K$. Hence the K\L-inequality (\ref{KL})  holds on a neighborhood $U$ of $K$. Since there are only finitely many 
$x_k$ outside $U$, on re-labeling the sequence, we may assume that (\ref{KL}) holds for the entire sequence.

By concavity of the de-singularizing function  $\phi$ in (\ref{KL}) we have
\begin{align}
\label{long}
\begin{split}
\phi(f(x_k)-f^*)  - \phi(f(x_{k+1})-f^*)
&\geq \phi'( f(x_{k})-f^*)( f(x_k)-f^* -(f(x_{k+1})-f^*)) \\
&= \phi'( f(x_{k})-f^*)( f(x_k) -f(x_{k+1})) \\
&\geq  \phi'( f(x_{k})-f^*) \lambda_k^{-1} \Psi(x_{k+1},x_k) \\
&\geq  \phi'( f(x_{k})-f^*) \lambda_k^{-1} m \|Px_{k+1}-Px_k\|^2.
\end{split}
\end{align}
Here the third line uses (\ref{descent2}), while the last line
uses the lower norm bound on $V$, which as mentioned above is based on Lemma \ref{lower}, applied to $\Psi_V$.

5)
By the Kurdyka-\L ojasiewicz inequality (\ref{KL})  we have
$$
\phi'(f(x_k)-f^*)  \|g_k\| \geq \gamma,
$$
using $g_k \in \partial f(x_k)$.  Applying the partial upper norm bound and  approximate optimality
at stage $k-1$, 
we get
\begin{align}
\label{also_take}
\begin{split}
\phi'(f(x_k)-f^*)^{-1} \leq \gamma^{-1} \|g_k\| &=
\gamma^{-1} \left( \|e_k\|+ \lambda_{k-1}^{-1} \|\nabla_1 \Psi(x_k,x_{k-1})\|\right)\\
&\leq \gamma^{-1}\left( \|e_k\| + M\lambda_{k-1}^{-1} \|Px_k-Px_{k-1}\|\right),
\end{split}
\end{align}
where  the upper norm bound on $V$ occurs in the last estimate. This bound uses Lemma \ref{upper} in the next section applied to $\Psi_V$.
Combining this with (\ref{long}) gives
$$
\phi(f(x_k)-f^*)  - \phi(f(x_{k+1})-f^*) \geq \frac{\gamma\lambda_k^{-1} m \|Px_{k+1}-Px_{k} \|^2}{\|e_k\| + M\lambda_{k-1}^{-1} \|Px_k-Px_{k-1}\|}.
$$
Since $a^2 \leq bc$ for $a,b,c \geq 0$ implies $a \leq \frac{1}{2}b + \frac{1}{2}c$, 
we get
\begin{align}
\label{26}
\begin{split}
\|Px_{k+1}-Px_k\| \leq &\textstyle\frac{1}{2} \left( \|Px_k-Px_{k-1}\| + M^{-1}\lambda_{k-1} \|e_k\| \right)
\\
&+\textstyle\frac{1}{2} \gamma^{-1} \lambda_{k}/\lambda_{k-1} M/m \left(\phi(f(x_k)-f^*)-\phi(f(x_{k+1})-f^* )\right)
\end{split}
\end{align}
and setting $C=\frac{Mr}{\gamma m}$ while using $\lambda_k/\lambda_{k-1} \leq r$, we get
\begin{equation}
\label{take}
\|Px_k-Px_{k+1}\| \leq \frac{1}{2}  \|Px_k-Px_{k-1}\| +\frac{\lambda_{k-1}}{2M}\|e_k\| + \frac{C}{2} \left[\phi(f(x_k)-f^*)  - \phi(f(x_{k+1})-f^*)  \right]. 
\end{equation}
Summing this from $k=1$ to $k=n$ gives
$$
\sum_{k=1}^n \|Px_{k+1}-Px_k\| \leq \frac{1}{2}\sum_{k=1}^n \|Px_k-Px_{k-1}\| + \frac{1}{2M} \sum_{k=1}^n   \lambda_{k-1}\|e_k\| + \frac{C}{2} \left[ \phi(f(x_1)-f^*)  - \phi(f(x_{n+1})-f^*)  \right].
$$
Hence
\begin{align*}
\sum_{k=1}^n \|Px_{k+1}-Px_k\| &\leq \|Px_1-Px_0\| + C  \left[ \phi(f(x_1)-f^*)  - \phi(f(x_{n+1})-f^*)  \right]\\
&\hspace{3cm}  - \|Px_{n+1}-Px_n\| + \frac{1}{2M} \sum_{k=1}^n \lambda_{k-1}\|e_k\|\\
&\leq \|x_1-x_0\| + C \phi(f(x_1)-f^*) + \frac{1}{2M} \sum_{k=1}^\infty \lambda_{k-1} \|e_k\|.
\end{align*}
Under condition a. the series on the right converges, hence 
$\sum_k \|Px_{k+1}-Px_k\|<\infty$, so that $Px_k$ is a Cauchy
sequence, which converges to some $Px^*\in V$, where $x^*$ is an accumulation point of the $x_k$. In fact,
$Px^*$ is the same for all accumulation points $x^*$ of the $x_k$.

6)
It remains to show that at least one accumulation point $x^*$ is critical under condition a.
Since $x_k \in G$, we have
$$
e_{k+1} \in \partial f(x_{k+1}) + \lambda_k^{-1} \nabla_1 \Psi(x_{k+1},x_k),
$$
so we can write
$\lambda_ke_{k+1} = \lambda_k g_{k+1} + v_{k+1}$ for $g_{k+1} \in \partial f(x_{k+1})$ and $v_{k+1} =\nabla_1 \Psi(x_{k+1},x_k)$. 
Now $\|\nabla_1 \Psi(x_{k+1},x_k)\|\leq M \|Px_{k+1}-Px_k\|$ and $\sum_k \|Px_{k+1}-Px_k\| < \infty$, hence
$v\in \ell_1$, and since also $\lambda\cdot e\in \ell_1$ by hypothesis a., we must have $\lambda \cdot g \in \ell_1$. By assumption, $\lambda \not \in \ell_1$, and this
means $g$ cannot be bounded away from $0$. In other words, $g_{k'} \to 0$ for at least a subsequence $k'$,
and then $0 \in \partial f({x}^*)$ from $x_{k'} \to {x}^*$, $g_{k'} \in \partial f(x_{k'})$, $\partial f$ being upper semi-continuous.

7)
If the stronger $\lambda_k \geq \eta > 0$ holds, then every accumulation point $x^*$ is critical, because in that case we must have $g\in \ell_1 \subset c_0$ for the entire sequence. 

8) It remains to discuss conditions b. and c. Under b. we can directly get rid of the term $\|e_k\|$
in the estimate (\ref{also_take}), and the same goes for condition c. The remainder of the proof is then simplified as we can work as if $e_k=0$.
\hfill $\square$
\end{proof}

\begin{remark}
1)
For $P=I$, $\Psi(x^+,x) = \frac{1}{2}\|x^+-x\|^2$, $\eta \leq \lambda_k \leq R$ and $e_k=0$  this was proved by Attouch and Bolte \cite{attouch_bolte}.
See also \cite[Thm. 24]{bolte2}.

\noindent
2)
When $P=I$ then $\sum_k \lambda_{k-1}\|e_k\| < \infty$ suffices for the limit point $x^*$ to be critical.

\noindent
3)
Constraints $x\in M$ are included directly by letting $f = f_0 + i_M$.
Inf-compactness of $f = f_0 + i_M$ on $G$ assures boundedness of the sequence $x_k$ and
is trivially satisfied if $M\subset G$ is closed bounded.

\noindent
4)
We can dispense with the hypothesis of separatingness of $\Psi_V$ if we assume $Px_k-Px_{k-1} \to 0$.

\noindent 
5) If we assume that upper and lower norm bounds still hold with the same $m,M$ as iterates approach the boundary $\partial G$, then
convergence holds also at the boundary. However, in the context of EM this is not a realistic assumption.
\end{remark}

\subsection{Lemmas for convergence}
Recall that the general regularizer $\Psi(x^+,x)$ is of class $C^2$ in the first variable
on $G\times G$, where $G \subset {\rm dom} \Psi(\cdot,y)$ for every $y\in G$, and
$G \subset {\rm dom}\Psi(x,\cdot)$ for every $x\in G$. Moreover, $\nabla_1\Psi(x,y)$ and
$\nabla^2_{11}\Psi(x,y)$ are jointly continuous. 
We have $\Psi(\cdot,x) \geq 0$ and $\Psi(x,x)=0$,
hence $\nabla_1 \Psi(x,x) = 0$ and $\nabla_{11}^2 \Psi(x,x) \succeq 0$ from the necessary optimality conditions.

\begin{lemma}
{\bf (Upper norm bound)}.
\label{upper}
    Let $K \subset G$ be compact convex. Then there exist $M > 0$ and $\delta > 0$ such that
    $\|\nabla_1 \Psi(x,y)\| \leq M\|x-y\|$ for all $x,y\in N(K,\delta)$.
\end{lemma}
\begin{proof}
Using compactness of $K\times K$  and continuity of $(x,y) \mapsto \nabla^2_{11}\Psi(x,y)$ on $G\times G$, choose $M>0$ such that $\lambda_{\max}(\nabla^2_{11} \Psi(x,y)) \leq M/2$ for all
$x,y\in K$. Then  find a neighborhood $N(K,\delta)$ of $K$ such that $\lambda_{\max}(\nabla^2_{11} \Psi(x,y)) \leq M$
for all
$x,y\in N(K,\delta)$.
Fix $\|h\|=1$, then Taylor expansion of $x \mapsto  h\cdot \nabla_1\Psi(x,y)$ at $y$ gives
$$
h\cdot \nabla_1 \Psi(x,y) =h\cdot \nabla_1\Psi(y,y) + h\cdot \nabla^2_{11}\Psi(\bar{y},y)(x-y)
$$
for some $\bar{y} \in (x,y)$, the open segment, where $\bar{y}$ depends on $x,y$ and $h$. 
Using $\nabla_1\Psi(y,y)=0$ and $\lambda_{\max}(\nabla^2_{11} \Psi(\bar{y},y)) \leq M$, we obtain
$h\cdot \nabla_1\Psi(x,y)  \leq \|h\| M \|x-y\|=M\|x-y\|$, which gives the claimed estimate, since $\|h\|=1$ is arbitrary.
\hfill $\square$
\end{proof}

We call this the upper norm bound. 
When $\nabla^2_{11} \Psi$ is strictly positive, we have the following

\begin{lemma}
\label{lower}
{\bf (Lower norm bound)}.
Let $K$ be a compact convex subset of $G$.  Suppose $\nabla^2_{11} \Psi(x,x) \succ 0$ on $K$.
Then there exist $\delta > 0$ and $m > 0$ such that
    $m\|x-y\|^2 \leq \Psi(x,y)$
for all $x,y$ with  $y\in N(K,\delta)$ and $\|x-y\| < \delta$.
\end{lemma}

\begin{proof}
1)
Since $K$ is compact and $\nabla^2_{11}\Psi(x,x) \succ 0$ on $G$, there exists $m>0$
such that $\nabla^2_{11}\Psi(x,x) \succeq m>0$ for all $x\in K$. From that we obtain a neighborhood $U=N(K,\delta)$ of $K$
such that $\nabla^2_{11}\Psi(y,z) \succeq m/2$ for all $y,z$ with $z\in N(K,\delta)$ and $\|y-z\| < \delta$. 

Indeed, for $x\in K$ choose $\epsilon_x > 0$ such that
$\nabla^2_{11} \Psi(y,z) \succeq m/2$ for all $y,z\in B(x,\epsilon_x)$. This is possible due to continuity of 
$(y,z) \mapsto \nabla^2_{11} \Psi(y,z)$. Now let $\Delta_K = \{(x,x): x\in K\}$
be the diagonal, then $\Delta_K \subset \bigcup_{x\in K} B(x,\epsilon_x/2)\times B(x,\epsilon_x/2)$, hence by
compactness of $\Delta_K$  there are finitely many $x_i \in K$ such that
$\Delta_K \subset B(x_1,\epsilon_1/2)\times B(x_1,\epsilon_1/2) \cup \dots \cup B(x_n,\epsilon_n/2) \times B(x_n,\epsilon_n/2)$,
where  $\epsilon_i = \epsilon_{x_i}$.

Now let $\delta := \min_{i=1,\dots,n}\epsilon_i/4$. Suppose
$z\in N(K,\delta)$ and $\|y-z\| < \delta$. Find $x\in K$ with $\|z-x\| < \delta$, then $\|y-x\| < 2 \delta$.
For some $x_i$ we have $(x,x) \in B(x_i,\epsilon_i/2) \times B(x_i,\epsilon_i/2)$, therefore
$\|y-x_i\| < 2\delta + \epsilon_i/2 < \epsilon_i$ and $\|y-x_i\| < \delta + \epsilon_i/2 < \epsilon_i$.
Hence $\nabla^2_{11} \Psi(y,z) \succeq m/2$ by the definition of $B(x_i,\epsilon_{x_i})$.

2)
Now  second order Taylor expansion of $\Psi(\cdot,y)$  at $y$ gives
$$
\Psi(x,y) = \Psi(y,y) + \nabla_1 \Psi(y,y)\cdot (x-y) + \textstyle\frac{1}{2} (x-y) \cdot \nabla^2_{11} \Psi(\bar{y},y)(x-y)
$$
for some $\bar{y} \in (x,y)$, the open segment. Therefore, if $\|x-y\| < \delta$, and $y\in N(K,\delta)$, 
then also $\|\bar{y}-y\| < \delta$, hence by part 1), 
$\lambda_{\min}(\nabla^2_{11}\Psi(\bar{y},y)) \geq m/2$, and
then using $\Psi(y,y)=0$ and $\nabla_1\Psi(y,y)=0$, we get
$\Psi(x,y) \geq (m/4) \|x-y\|^2$.
\hfill $\square$
\end{proof}

With the proof of Lemma \ref{lower} we can also get the following

\begin{lemma}
\label{lem9}
 Let $K \subset G$ be compact and let
    $\Psi(x^+,x)$ have a lower norm bound at every $x\in K$. Then there exist $\delta > 0$ and $m$  such that
$m\|y-z\|^2 \leq \Psi(y,z)$ holds for all $y,z\in N(K,\delta)$ with $\|y-z\| < \delta$. 
\end{lemma}

\subsection{Rate of convergence}
\label{sect_rate}
\begin{corollary}
\label{complex}
Consider the case $\lambda_k \geq \eta > 0$ and $e_k=0$ in Theorem {\rm \ref{main}}.
Further suppose that
  $\phi(s) = s^{1-\theta}/(1-\theta)$   
  for $\theta \in [\frac{1}{2},1)$.
  If $\theta\in (\frac{1}{2},1)$, then
  the speed of convergence is $\|Px_k-Px^*\| = O(k^{-\frac{1-\theta}{2\theta-1}})$. For 
  $\theta=\frac{1}{2}$ the speed is R-linear.
  \end{corollary}
  
  \begin{proof}
In the \L ojasiewicz case equation (\ref{take}) specializes to
$$
\|Px_k-Px_{k+1}\| \leq \frac{1}{2} \|Px_k-Px_{k-1}\| + \frac{C}{2} \left[(f(x_k)-f^*)^{1-\theta}  - (f(x_{k+1})-f^*)^{1-\theta}  \right]. 
$$
Summing this form $k=N$ to $k=M$ gives
\begin{align*}
&-\frac{1}{2} \|Px_{N-1}-Px_N\| + \frac{1}{2} \sum_{k=N}^{M-1}\|Px_k-Px_{k+1}\| + \|Px_M-Px_{M+1}\|    \\
&\hspace{5cm}\leq \frac{C}{2} \left[(f(x_N)-f^*)^{1-\theta}  - (f(x_{M+1})-f^*)^{1-\theta}  \right].
\end{align*}
Passing to the limit $M\to \infty$ gives
\begin{align*}
-\frac{1}{2} \|Px_{N-1}-Px_N\| + \frac{1}{2} \sum_{k=N}^{\infty}\|Px_k-Px_{k+1}\| 
\leq \frac{C}{2} (f(x_N)-f^*)^{1-\theta}.
\end{align*}
Putting $S_N = \sum_{k=N}^\infty\|Px_{k} - Px_{k+1}\|$, this becomes
$$
-\frac{1}{2}(S_{N-1}-S_N) + \frac{1}{2} S_N \leq \frac{C}{2} (f(x_N)-f^*)^{1-\theta}.
$$
Now from (\ref{also_take}) we have $\phi'(f(x_N)-f^*)^{-1} \leq \gamma^{-1} M\lambda_{k-1}^{-1} \|Px_N-Px_{N-1}\| =
\gamma^{-1}  M\lambda_{k-1}^{-1}(S_{N-1}-S_N) \leq \gamma^{-1}M\eta^{-1} (S_{N-1}-S_N)$.
Since $\phi'(s)=s^{-\theta}$, this implies
\begin{align*}
\phi(f(x_N)-f^*)&=(1-\theta)^{-1}(f(x_N)-f^*)^{1-\theta} \\&= (1-\theta)^{-1}\left[ \phi'(f(x_N)-f^*)^{-1}\right]^\frac{1-\theta}{\theta}\\
&\leq (1-\theta)^{-1}  (M\gamma^{-1}\eta^{-1})^\frac{1-\theta}{\theta}(S_{N-1}-S_N)^\frac{1-\theta}{\theta}.
\end{align*}
So altogether we get
\begin{equation}
\label{rpm}
\frac{1}{2}S_N \leq C'(S_{N-1}-S_N)^\frac{1-\theta}{\theta} + \frac{1}{2}(S_{N-1}-S_N)
\end{equation}
for $C'=(1-\theta)^{-1}  (M\gamma^{-1}\eta^{-1})^\frac{1-\theta}{\theta}$. Now for $\theta > \frac{1}{2}$ we have $\frac{1-\theta}{\theta} < 1$, so the first term on the right 
of (\ref{rpm}) dominates the second term, and we get
$$
S_N^\frac{\theta}{1-\theta} \leq C''(S_{N-1}-S_N)
$$
for $N$ large enough and yet another constant $C''$. Following \cite[Cor. 4(24)ff]{noll}, this leads to an estimate
$S_N \leq C''' N^{-\frac{1-\theta}{2\theta-1}}$. 

It remains to discuss the case $\theta=\frac{1}{2}$. Here (\ref{rpm}) gives
$$
\frac{1}{2}S_N \leq C'(S_{N-1}-S_N) + \frac{1}{2}(S_{N-1}-S_N),
$$
hence
$$
S_N \leq \frac{1+SC'}{2+2C'} S_{N-1}
$$
which gives Q-linear convergence of the $S_N$, hence R-linear convergence of the $Px_k$. 
\hfill $\square$
\end{proof}

\begin{remark}
When $\phi(s) = s^{1/2}$, which is the best possible case, the \L ojasiewicz inequality (\ref{KL}) specializes to the Polyak-\L ojasiewicz inequality. 
In the unconstrained case, for a critical point $\bar{x}$,  (\ref{KL}) is then $(f(x) - f(\bar{x}) )^{-1/2} \|\nabla f(x)\| \geq \gamma$, and that
means $f$ is locally bounded below by a quadratic, and in particular, has a strict local minimum at $\bar{x}$. Indeed, assuming $f(\bar{x})=0$, $\bar{x}=0$
and letting $y(t) = f(th)$ for fixed $\|h\|=1$, (\ref{KL}) gives $y' \geq \gamma \sqrt{y}$, hence $dy/\sqrt{y} \geq \gamma dt$, hence
$\sqrt{y} \geq \frac{\gamma}{2} t$, i.e. $y \geq \frac{\gamma^2}{4} t^2$, using $y(0)=0$. This is of course too good to be true, so
we expect the Polyak-\L ojasiewicz inequality to 
be satisfied in exceptional cases only.
\end{remark}

In the general case of the Kurdyka-condition we can still say something:
\begin{corollary}
\label{kurd}
Consider the case $\lambda_k\geq \eta > 0$ and $e_k=0$ in Theorem {\rm \ref{main}}. Then the speed of partial convergence is
$\| Px_k-Px^*\| \leq C\phi( f(x_k)-f^*)$ with $\phi$ the de-singularizing function.
\end{corollary}

\begin{proof}
Since $\Psi(x,x)=0$ and $\nabla_1\Psi(x,x)=0$, 
Taylor expansion of $\Psi(\cdot,x)$  at $x$ gives $\Psi(u,x) = \Psi(x,x) + \nabla_1 \Psi(x,x)\cdot (u-x) + \frac{1}{2} (u-x) \cdot \nabla^2_{11} \Psi(\bar{x},x)(u-x)
= \frac{1}{2} (u-x) \cdot \nabla^2_{11} \Psi(\bar{x},x) (u-x)$
for $\bar{x}$ on the segment $(x,u)$. Since for a compact set  $K\subset G$ we find 
$c > 0$ such that $\nabla^2_{11} \Psi(x',x) \preceq c^2I$ for all $x',x\in K$, choosing as $K$ the convex hull of the set of iterates
and its accumulation points, we have the estimate $\Psi(x^+,x)^{1/2} \leq c \|x^+-x\|$ for the sequence generated by (\ref{general}).

Since the solution  $x^+$ of (\ref{general}) from the current $x$ is exact, letting
$f(x^+) = r^+$, $f(x) = r$,  we have
$\Psi(x^+,x) \leq \Psi(u,x)$ for all $u\in \{f \leq r^+\}$. Hence
$\Psi(x^+,x)^{1/2} \leq \min_{u\in \{f\leq r^+\}} \Psi(u,x)^{1/2} \leq  c \min_{u\in \{f\leq r^+\}}\|u-x\| = c\, {\rm dist}(\{f\leq r^+\},x)
\leq c \max_{v\in \{f\leq r\}} {\rm dist}(\{f\leq r^+\},v) \leq c \,{\rm haus}(\{ f\leq r^+\},\{f\leq r\})$, where we used that $x$ is one of the $v\in \{f\leq r\}$.

According to \cite[Thm. 20(vi),Cor. 4]{bolte2} the K\L-property is equivalent to Lipschitz continuity of the sublevel operator. More precisely,
there exists $k>0$ such that
${\rm haus}(\{f\leq r^+\},\{f\leq r\}) \leq k |\phi(r^+)-\phi(r)|$, where $\phi$ is the de-singularizing function of (\ref{KL}).
Substituting this gives
$\Psi(x^+,x)^{1/2} \leq c \cdot{\rm haus}(\{f\leq r^+\},\{f\leq r\}) \leq ck |\phi(r^+)-\phi(r)|$.
Using the lower norm bound, we deduce
$\|Px_{k+1}-Px_k\| \leq m^{-1/2}ck [ \phi(f(x_{k})-f^*) - \phi(f(x_{k+1})-f^*)]$, and summing both sides from $k=n$ to $k=r$
gives $\sum_{k=n}^r \|Px_{k+1}-Px_k\| \leq m^{-1/2}ck [\phi(f(x_n)-f^*) - \phi(f(x_{r+1})-f^*)|$. Letting $r\to \infty$
gives the claimed rate
$\|Px_n-Px^*\| \leq m^{-1/2}ck \phi(f(x_n)-f^*)$.
\hfill $\square$
\end{proof}

For the Kurdyka inequality we may also argue as follows.
From (\ref{take}) we  obtain the estimate
$$
-\textstyle\frac{1}{2}(S_{n-1}-S_n) + \frac{1}{2}S_n \leq \frac{C}{2} \phi( f(x_n)-f^*).
$$
Fix $\alpha \in (\frac{1}{2},1)$, and divide integers in two classes
$\mathcal N_1 = \{n: S_n \leq \alpha S_{n-1}\}$ and $\mathcal N_2 = \{n: S_n > \alpha S_{n-1}\}$.
Now for $n \in \mathcal N_2$ we have
$$
(1-\textstyle\frac{1}{2\alpha}) S_n \leq S_n - \frac{1}{2}S_{n-1} \leq \frac{C}{2} \phi(f(x_n)-f^*).
$$
On the other hand, for $n\in \mathcal N_1$ we have $S_n \leq \alpha S_{n-1}$, so here the error shrinks
with linear rate.  Altogether we get a sequence $n_1 < m_1 < n_2 < m_2 < \dots$ such that
\begin{align*}
S_n &\leq c \phi(f(x_n)-f^*) \mbox{ for $n=n_k,\dots,m_k-1$} \\ 
S_{n} &\leq \alpha^{n-m_k} S_{m_k},\quad\;\,\, \mbox{ for $n = m_k,\dots,n_{k+1}-1$},
\end{align*}
with $c=\frac{C}{2}(1-\frac{1}{2\alpha})^{-1}$, 
which is a slight refinement of Corollary \ref{kurd}, as it leaves the option of the entire sequence $S_n$ converging
linearly with rate $\alpha$ even  for $\phi$ less desingularizing that $s^{1/2}$.

\section{Convergence of the EM algorithm}
\label{convergence_EM}
Now we apply this to the EM algorithm. In the first place, we assume interiority.

\begin{theorem}
\label{theorem2}
{\bf (Convergence for constrained exponential family)}.
Let $q(y,\theta)$ be incomplete data from a minimal $n$-dimensional exponential family $p(x,\theta)$.
Suppose $-\log q(y,\cdot)$ and $M$ are definable.     Let $\theta^{(k)}\in M$ be a bounded sequence generated by the
constrained  {\rm EM} algorithm which together with its accumulation points stays in $G$.
    Then the sequence $\theta'^{(k)}$ of accurate parameters converges, $\theta'^{(k)} \to \theta'^*$. Every accumulation point
    $\theta^*$ of the sequence $\theta^{(k)}$ solves the constrained incomplete data MLE problem,  has the same projection $P\theta^*=\theta'^*$, and the    
    conditional distributions $\mathbb P_{\theta^{(k)}}^{x|y}$ converge weakly to   $\mathbb P_{\theta^*}^{x|y}$, the limit being the same for all $\theta^*$.
   Moreover, the sequence $\xoverline{T}(x_k)$ is also convergent. 
           \end{theorem}

\begin{proof}
    Since $-\log q(y,\cdot)$ and $M$ are definable, so is the  objective $f = -\log q(y,\cdot) + i_M$ in (\ref{missing}). Hence $f$ has the K\L-property (\ref{KL}). 
 By hypothesis the sequence together with its accumulation points stays in $G$,
 as required for our main convergence result.

    By Proposition \ref{prop4} we have $\nabla^2_{\theta'\theta'}\bar{\psi}_y(\theta') \succ 0$ for the accurate parameter $\theta'$, hence
    the partial regularizer $\Psi(\theta^+,\theta) = \bar{K}_y(\theta'||\theta'^+)$ satisfies $\nabla^2_{22} \bar{K}_y(\theta'||\theta') \succ 0$. Therefore
    it has a lower norm bound on the compact set of accumulation points of the $\theta'^{(k)}$ contained in the subspace $V$
    of dimension $m$. Also, since the family $\bar{k}(\cdot|y,\theta')$ is minimal,
    the partial regularizer is separating on the subspace $V$.
    Therefore we can apply the main convergence theorem (with $\lambda_k=1$), and this gives convergence of the $\theta'^{(k)}$.
    
    Since $\mathbb P_{\theta^{(k)}}^{x|y} = \mathbb P_{\theta'^{(k)}}^{x|y}$ by Proposition \ref{prop4}, and since the right hand sequence converges
    weakly, so does the left hand sequence. It also follows from a classical result of F. Riesz  that $\bar{k}(\cdot|y,\theta'^{(k)})$ converges to $\bar{k}(\cdot|y,\theta'^*)$ in 
    $L^1(h^{-1}(y),d\mu_y)$. But then due to (\ref{20}) the sequence $k(\cdot|y,\theta^{(k)})$  also converges in $L^1$, regardless of whether the $\theta''^{(k)}$
    converge.
    
    As a consequence of the main convergence theorem, every accumulation point $\theta^*$ of the full sequence $\theta^{(k)}$ is critical,
    and $P\theta^* = \theta'^*$ for every such $\theta^*$.
    Convergence of $\xoverline{T}(x_{k+1})$ follows because $\xoverline{T}(x_{k+1})$ is the result of computing
    $\mathbb E_{\theta'^{(k)}}[\xoverline{T}(x)|y]$, which depends continuously on $\theta'^{(k)}$.  
        \hfill $\square$
\end{proof}

\begin{remark}
Partial
convergence is also guaranteed under the more general choices of $\lambda_k$ in Theorem \ref{main}, but that
requires computing $-\log q(y,\theta^+)$ and $\bar{K}_y(\theta'||\theta'^+)$ separately. In practice one prefers to  use
 the Q-function. 
\end{remark}

\begin{corollary}
\label{cor3}
Under the hypotheses of Theorem {\rm \ref{theorem2}}, suppose $T(h^{-1}(y))$ is not contained in a proper affine subspace of $\mathbb R^n$.
Then the entire sequence $\theta^{(k)}$ converges.
\end{corollary}

\begin{proof}
In this situation no dimension reduction takes place and we have $\theta^{(k)} = \theta'^{(k)}$. Hence $\theta^{(k)}$ converges.
\hfill $\square$ 
\end{proof}

\subsection{Trouble at the boundary}

From (\ref{q}) we have dom$\,q(y,\cdot) = \Theta$, while (\ref{k}) gives
$\Theta \subset {\rm dom} \Psi(\cdot,\theta')$ for every $\theta'\in G$. 
But this does not exclude the possibility that
iterates lie on the boundaries
of dom$(-\log q(y,\cdot))$ and dom$\,\Psi(\cdot,\theta')$ simultaneously, causing problems (see Example \ref{ex_missing}).  We take a closer look.

\begin{definition}
{\rm 
We call  the incomplete data problem {\it regular} if  $f(\theta)=-\log q(y,\theta)+i_M(\theta) =\infty$
for all $\theta \in \partial \Theta$, and we call it {\it steep} if
$\partial f(\theta)=\partial \left( -\log q(y,\cdot)+i_M\right)(\theta)= \emptyset$ for all $\theta\in \partial\Theta$. 
}
\end{definition}

\begin{remark}
1) We have $\partial f(\theta)=\partial\left(-\log q(y,\cdot)+i_M   \right)(\theta)= \partial (-\log q(y,\cdot))(\theta) + \mathcal N_M(\theta)$ for $\theta\in G$
by \cite[8.8. c]{rock} or \cite{mord}, but the sum rule fails at the boundary $\partial \Theta$, unless additional regularity hypotheses are made. This is why
the definition uses $\partial\left(-\log q(y,\cdot)+i_M   \right)$.

2) Regular implies steep. If $M \subset G = {\rm int} (\Theta)$,
then the problem is automatically regular. 

3) When $q(y,\cdot)$
is from an exponential family and $M=\xoverline{\Theta}$, then steepness  in the sense of the definition is equivalent to steepness of
the  log-normalizer $\psi_q$ of $q$.  

4) Suppose the complete data family is steep in the sense that whenever $\theta\in \partial\Theta$, then
$\partial_\theta \left( p(x,\cdot)+i_M\right) (\theta)= \emptyset$ for 
$\mu_y$-a.a. $x\in h^{-1}(y)$. Then steepness of the incomplete data problem follows
from the inclusion
$\partial_\theta \left( q(y,\cdot)+i_M\right)(\theta)\subset \int_{h^{-1}(y)} \partial_\theta \left(p(x,\cdot)+i_M\right)(\theta) d\mu_y(x)$
(see \cite[Thm. 2.7.2]{clarke}).  This happens when $p(x,\theta)$  is of exponential type, as then
$\partial_\theta p(x,\theta) = p(x,\theta) (T(x) + \partial \psi(\theta))$, so that 
$\partial \psi(\theta)=\emptyset$ for $\theta\in \partial\Theta\cap M$ forces
$\partial_\theta p(x,\theta)=\emptyset$ for all $x$.
This justifies the definition.
\end{remark}

Let $\theta^{(k)}$ be the sequence generated by the constrained EM algorithm,
then every accumulation point $\theta^*$ must have finite value $f(\theta^*) < \infty$, because
$\theta^{(k)} \in \{ f\leq f(\theta^{(1)})\}$, and by lower semi-continuity
of $f$ this set is closed. Therefore in the regular case no $\theta^*$
can be on the boundary $\partial\Theta$. Hence the interiority hypothesis in Theorem
\ref{theorem2} is automatically satisfied for bounded
$\theta^{(k)}$. Boundedness
is  assured e.g. under inf-compactness of $f$. In other words, if the incomplete data problem is regular,
there ain't any trouble at the boundary $\partial \Theta$.

Now
suppose $f$ is steep, and let  $e_k=0$.
By optimality $-\lambda_{k-1}^{-1}\nabla_1\Psi(\theta^{(k)},\theta^{(k-1)})\in \partial f(\theta^{(k)})$,
hence $\theta^{(k)}$ cannot be on the boundary 
$\partial \Theta$, and the method is well-defined. 
Suppose $\theta^{(k)} \to \theta^*\in \partial \Theta$  for a subsequence $k\in N\subset \mathbb N$. 
Then by steepness $\theta^*$ cannot be a solution of the constrained incomplete
data MLE program, because that would give $0 \in \partial f(\theta^*)
= \partial (-\log q(y,\cdot) + i_M)(\theta^*)=\emptyset$. 
Unfortunately, steepness alone does {\it not} prevent iterates from
approaching $\partial \Theta$, but at least we know in that case that these iterates
go astray.

\begin{remark}
When $M \subset G = {\rm int}(\Theta)$ is bounded, then
the sequence
$\theta^{(k)}$ together with its accumulation points stays in $G$ and all
trouble at the boundary $\partial \Theta$ is avoided. Conversely, suppose we
wish to make a statement
about a sequences $\theta^{(k)}$ respecting interiority, i.e., bounded and contained in $G$
together with its accumulation points. Then
we can replace the constraint set $M\subset \xoverline{\Theta}$
by a closed subset $M' \subset M \cap G$ so that the sequence may be considered
as generated under the constraints $\theta^{(k)} \in M'$. This set $M'$ can be chosen bounded,
definable if $M$ is definable, and convex because $G$ is convex. See
\cite[Sect. 2.2]{bregman}.
\end{remark}

\subsection{Consequences for the M step}
\label{sect_M}
Having established partial convergence $\theta'^{(k)} \to \theta'^*$ under interiority, we now investigate under what conditions we may
upgrade this to
convergence of the full parameter sequence $\theta^{(k)}$.
Observe that in the new coordinates $Q\theta = (\theta',\theta'')$ the M-step (\ref{missing}) has the equivalent form:
$$
({\theta'}^{(k+1)},{\theta''}^{(k+1)}) \in \argmin_{{\theta'},{\theta''}} -\log q(y,{\theta'},{\theta''}) + i_M({\theta'},{\theta''}) + \lambda_k^{-1}\Psi({\theta'},{\theta'}^{(k)}), $$
with
$$
\Psi({\theta'},{\theta'}^{(k)}) = -\mathbb E_{\theta'^{(k)}}\left[ \log \frac{\bar{k}(x|y,{\theta'})}{\bar{k}(x|y,{\theta'}^{(k)})} \bigg| y \right]
= \bar{K}_y({\theta'}^{(k)} || {\theta'}),
$$
a partial regularizer independent of the variable $\theta''$. Identifying for simplicity $\theta$ with $(\theta',\theta'')$, 
define $M(\theta') = \{\theta'': (\theta',\theta'')\in M\}$, then we can
split the M step optimization program as follows:
$$
(P_k) \quad \theta''^{(k+1)} \in \argmin_{\theta'' \in M(\theta'^{(k+1)})} -\log q(y,\theta'^{(k+1)},\theta''),
$$
which is just a sequence of parametrized optimization programs in $\theta''$, with $\theta'^{(k+1)}$ the parameter,  and no longer any regularization
affecting $\theta''$. The limiting
program is clearly:
$$
(P_\infty) \quad \theta''^* \in \argmin_{\theta'' \in M(\theta'^*)} -\log q(y,\theta'^*,\theta''),
$$
and since $\theta'^{(k)} \to \theta'^*$, every accumulation point $(\theta'^*,\theta''^*)$ of the EM sequence $\theta^{(k)}$ gives
a solution $\theta''^*$ of $(P_\infty)$ with the same incomplete data MLE value $q^* = q(y,\theta^*)$.
This has the following immediate consequence:

\begin{proposition}
If the limiting program $(P_\infty)$ has a unique critical point $\theta''^*\in M(\theta'^*)$ among those with critical value $q^*$, then the {\rm EM} sequence $\theta^{(k)}$ converges.
\end{proposition}

This is a weaker hypothesis than requesting as e.g. in \cite{wu}, that the full incomplete data MLE program
has a unique solution $\theta^*$ with the correct MLE value.

\begin{remark}
Note, however,  that in each program $(P_k)$ we are free to choose {\it any} of the local solutions in case there are several. If there exists $\epsilon > 0$ such that in every
$(P_k)$ one can choose two local solutions a distance $\epsilon$ apart, then failure of convergence of the spare sequence $\theta''^{(k)}$ is inevitable. This happens for instance in the counterexample in
\cite[Sect. 4]{vaida}. \end{remark}

A second consequence is based on the following.

\begin{proposition}
\label{principle}
Suppose the complete data exponential family $p(x,\theta)$ is minimal.
Then
$-\nabla^2_{{\theta''}{\theta''}} \log q(y,{\theta'},{\theta''}) \succ 0$ for fixed ${\theta'}$.
\end{proposition}

\begin{proof}
Using standard notation, one defines
$I(\theta,y) = -\nabla^2_{\theta\theta} \log q(y,\theta)$, which makes $\mathcal I(\theta) = \mathbb E_\theta[I(\theta,y)]$ the expected Fisher information of  incomplete data. 
For complete data one defines
$I_c(\theta,x) = -\nabla^2_{\theta\theta} \log p(x,\theta)$, then 
$\mathcal I_c(\theta) = \mathbb E_\theta[I_c(x,\theta)]$ is the expected Fisher information of complete data. 
In the same vein, one also lets 
$\mathcal I_c(\theta,y) = \mathbb E_\theta( I_c(\theta,x)| y)=\nabla^2 \psi(\theta)$,  the conditional expected Fisher information of complete data given $y$.

Now from (\ref{k})
we get $\log p(x,\theta) = \log q(y,\theta) + \log k(x|y,\theta)$, hence differentiating twice gives
$$
I_c(\theta,x) = I(\theta,y) - \nabla^2_{\theta\theta} \log k(x|y,\theta).
$$
Taking conditional expectations over $x$ given $y$, we obtain
\begin{equation}
\label{add}
\mathcal I_c(\theta,y) = I(\theta,y) + \mathcal I_m(\theta,y),
\end{equation}
where
$\mathcal I_m (\theta,y) = \mathbb E_\theta(-\nabla^2_{\theta\theta} \log k(x|y,\theta) | y)$ is the expected Fisher information of missing data 
conditioned on $y$.
The latter, however, was previously identified as the second derivative of the regularizer.

Now for $k(x|y,\theta)$
an exponential family, and adopting the change of coordinates in (\ref{Q}), we know that 
$$
\mathcal I_m(\theta,y) = \left[\begin{array}{c|c} \nabla^2_{\theta'\theta'} \bar{\psi}_y(\theta') & 0  \\ \hline 0&0\end{array}\right], \qquad \nabla^2_{\theta'\theta'} \bar{\psi}_y(\theta') \succ 0,
$$ 
while
$$
I(\theta,y)=
-\nabla^2_{\theta\theta} \log q(y,\theta',\theta'')=\left[
\begin{array}{c|c} -\nabla^2_{\theta'\theta'}  \log q(y,\theta',\theta'') &  -\nabla^2_{\theta'\theta''}  \log q(y,\theta',\theta'') \\ \hline
*&  -\nabla^2_{\theta''\theta''}  \log q(y,\theta',\theta'') \end{array}\right].
$$
Since $p(x,\theta)$ is minimal by hypothesis, we have $I_c(\theta,x) = \nabla^2 \psi(\theta) \succ 0$.
Therefore also $\mathcal I_c(\theta,y) \succ 0$.
In consequence,  the matrices $\mathcal I_m(\theta,y)$ and $I(\theta,y)$ on the right of (\ref{add}) must add up to a matrix of full rank.  Due to the structure of $\mathcal I_m(\theta,y)$
this forces $-\nabla^2_{\theta''\theta''} \log q(y,\theta) \succ 0$ for the lower diagonal block in $I(\theta,y)$. 
\hfill $\square$
\end{proof}

Applying this with $\theta'^*$ shows that the objective function
$-\log q(y,\theta'^*,\cdot)$ of $(P_\infty)$ is strictly convex. We therefore have the following consequence:

\begin{theorem}
\label{theorem3}
Under the assumptions of Theorem {\rm \ref{theorem2}}, suppose $M(\theta'^*)$ is convex. Then the EM sequence $\theta^{(k)}$ converges. 
\end{theorem}

\begin{proof}
A strictly convex function has a unique minimum on a convex domain.
\hfill $\square$
\end{proof}

Note that $M(\theta'^*)$ is clearly convex if $M$ is convex, but convexity of $M(\theta'^*)$ is a weaker hypothesis. In particular
we have the following

\begin{corollary}
\label{cor4}
Suppose the constraint set is
 $M=\{\theta(u): u \in U\} = \{(\theta'(u),\theta''(u)): u \in U\}$, with $U \subset \mathbb R^m$, $\theta'(\cdot)$ of class $C^1$
and $\theta''(\cdot)$ continuous, both definable. Let $\theta'(u^*)=\theta'^*$, and suppose
the rank of the Jacobian $\frac{d\theta'}{du}(u^*)$ is $m$. Then the {\rm EM} sequence $\theta^{(k)}$ converges. When $\theta''(\cdot)$ is locally Lipschitz, then
the speed of convergence of $\theta^{(k)}$ is the same as that of $\theta'^{(k)}$.
\end{corollary}

\begin{proof}
Under the rank hypothesis the mapping $u \mapsto \theta'(u)$ has locally a left inverse, i.e., we have a $C^1$
mapping $\theta' \mapsto u(\theta')$ defined in a neighborhood of $\theta'^*$ such that $u(\theta'^*)=u^*$ and
$u(\theta'(u))=u$. Then
$\theta'' = \theta''(u) = \theta''(u(\theta'))$, so that $\theta''$ is a function of $\theta'$. Then $M(\theta'^*)$ is singleton, hence convex.
An even more direct argument is that $\theta''^{(k)} = \theta''(u(\theta'^{(k)}))$ converges by continuity of $\theta''(\cdot)$ and $u(\cdot)$.
The statement concerning speed follows because when $\theta''(\cdot)$ is locally Lipschitz, then so is
$\theta''(\cdot) \circ u(\cdot)$.
\hfill $\square$
\end{proof}

\begin{remark}
1)
Assuming that each $(P_k)$ has a unique solution is not sufficient for convergence, as the $\theta''^{(k)}$ obtained may still have several accumulation points.

2)
Under the somewhat artificial assumption that the set of accumulation points of the sequence $\theta''^{(k)}$ is discrete, one obtains convergence
as soon as $\theta''^{(k)} - \theta''^{(k-1)} \to 0$.

3)
When dependence of the solution set $\argmin_{\theta''} -\log q(y,\theta',\cdot) + i_{M(\theta')}$ on the parameter $\theta'$ is upper Lipschitz (cf. \cite{klatte})
on the compact set $\{\theta'^{(k)}: k\in \mathbb N\} \cup \{\theta'^*\}$, then convergence follows, because this forces $\|\theta''^{(k)} - \theta''^{(k-1)}\|
\leq L \|\theta'^{(k)} - \theta'^{(k-1)}\|$ for some $L>0$, and since $\sum_k \|\theta'^{(k)}-\theta'^{(k-1)}\| < \infty$, the sequence $\theta''^{(k)}$ is also Cauchy.
For NLP constraints $M$, sufficient conditions are  discussed in \cite{robinson,klatte}, are 
typically local, and require mild regularity hypotheses. Here these have to be satisfied at all
accumulation points $\theta''^*$ of the spare sequence $\theta''^{(k)}$.
\end{remark}

\begin{theorem}
\label{theorem4}
{\bf (Regularized EM for exponential family)}. Under the hypotheses of Theorem {\rm \ref{theorem2}}, suppose the {\rm M} step
is regularized as 
$\min_{\theta \in M} - Q(\theta,\theta^{(k)}) + \lambda_k^{-1} \|\theta'' - \theta''^{(k)}\|^2$.
Then the sequence $\theta^{(k)}$ converges to a critical point $\theta^*$ which is a MLE for the incomplete data problem.
The value of the incomplete data negative log-likelihood  is still monotonically decreasing.
\end{theorem}

\begin{proof}
In view of Proposition \ref{prop1} we have modified the M step such that
the regularizer is now
$\Psi(\theta,\theta^{(k)})=\bar{K}_y(\theta'^{(k)}|| \theta') + \|\theta''-\theta''^{(k)}\|^2$, which is no longer partial but full. We apply the main convergence theorem with 
$P=I$, which gives convergence of the $\theta^{(k)}$. Since (\ref{general}) gives always decrease of the objective, the last  statement follows.
\hfill $\square$
\end{proof}

\subsection{Definable objectives}
We inquire whether, or when, objectives $f(\theta) = -\log q(y,\theta) + i_M(\theta)$ in (\ref{missing})  are definable in an o-minimal structure, because this is
how we assure
 the K\L-inequality (\ref{KL}). Starting with $-\log q(y,\theta)$, we first
run through those cases where $q(y,\theta)$ is by itself from an exponential family, say
$q(y,\theta) = \exp\{ \langle \theta,T(y) \rangle - \psi_q(\theta) + h(y)\}$. Here definability hinges on definability of the corresponding
log-normalizer $\psi_q$.

Inspecting lists of exponential families, one finds that along with algebraic expressions,  log-normalizers $\psi(\theta)$ sometimes include
terms of the form $\log \theta_i$ for certain components $\theta_i$ of $\theta$.
Those are definable in $\mathbb R_{an,exp}$. Moreover, if these $\theta_i$ can a priori be bounded
and bounded away from 0, one gets definability in $\mathbb R_{an}$. Inverse gamma distribution and $\chi^2$-distribution call for 
terms of the form $\log \Gamma(\theta_i)$ for certain components of $\theta$, which require 
definability
 of the Gamma function. This has recently been addressed e.g. in \cite{padgett}, and for our purpose it is again sufficient to
 bound these $\theta_i$ away from 0.

 The second case is when $q(y,\theta)$ are incomplete data from an exponential family, as termed in \cite{sundberg}, but do not by themselves
stem from an exponential family. Here  due to $-\log q(y,\theta) = -\log p(x_k,\theta) + \lambda_k^{-1} \bar{K}_y(\theta'^{(k)} || \theta')$
 definability of $\log q(y,\cdot)$ may be derived from definability of $\log p(x_k,\cdot)$
 in tandem with definability of $\bar{K}_y(\theta'^{(k)}|| \cdot)$. The first is assured when $\psi$ is definable, as discussed above. For the second
 we use Proposition \ref{Kull_Breg}, which shows that definability
 of the Bregman distance induced by $\psi_y$ is required, and this follows from definablity of $\psi_y$.

 Definability of $M$ is even less complicated, as $M$ typically
 gathers equality and inequality constraints of the form
 $M = \{\theta\in \mathbb R^n: f_i(\theta) = 0, i\in I, g_j(\theta) \leq 0, j\in J\}$ for finite sets $I,J$ and definable functions $f_i, g_j$, typically sub-analytic or even algebraic.
 Note that $M$ may even have the benefit to restrict components $\theta_i$ to a bounded interval, which
 allows to replace $\mathbb R_{an,exp}$ by the more convenient structure $\mathbb R_{an}$, where (\ref{KL}) turns into the \L ojasiewicz inequality.

\begin{remark}
When $f$ is definable in $\mathbb R_{an}$,  Corollary \ref{complex} gives a convergence
rate $\| \theta'^{(k)} - \theta'^*\| = O(k^{-\rho})$. If in addition $\theta''=\theta''(\theta')$ is locally Lipschitz,
we get the same rate for the full parameter sequence. This holds in
Corollary \ref{cor3} and Theorem \ref{theorem4}, but also in the case in Corollary \ref{cor4}. 
In Theorem \ref{theorem3}
it also holds due to $\nabla^2_{\theta'' \theta''} -\log q(y,\theta',\cdot) \succeq  \epsilon > 0$ for $\theta'$  in
the compact set $\{\theta'^*\} \cup \{\theta'^{(k)}: k\in \mathbb N\}$, 
provided $M$ is given by sufficiently smooth definable equality and inequality constraints, where the MFCQ is satisfied, see \cite{klatte}. 

Linear speed is obtained in the case $\phi(s) = s^{1/2}$, which corresponds to the Polyak-\L ojasiewicz inequality.
Unfortunately this is a rather strong hypothesis (see also Section \ref{sect_alter} for that aspect).
\end{remark}

\subsection{Reasons for failure}

We can now list the following reasons why the EM algorithm
for incomplete data from a constrained exponential family may fail to converge to critical points:
\begin{enumerate}
\item Iterates may be unbounded or tend to the boundary of $\Theta$. Once those are ruled out:
\item Convergence may still fail because $-\log q(y,\cdot)+i_M$ does not have the K\L-property. 
But even $f$ does have the K\L-property:
\item It may still happen that only the accurate parameter $\theta'^{(k)}$ converges, while the spare sequence $\theta''^{(k)}$ fails to converge.  This may
happen  if $M(\theta'^*)$ is not convex. 
\item 
But even when $M(\theta'^*)$ is convex,
including the unconstrained case, convergence of the $\theta''^{(k)}$  may still fail
because $-\log q(y,\theta'^*,\cdot)$ is not strictly convex. This may be the case because  $p(x,\theta)$ is not minimal. The latter may be avoided when setting up the problem.
\end{enumerate}

In curved families $M = \{ \theta(u): u \in U\}$, chances of convergence are paradoxically even better, as some of the degrees of freedom
are removed. As we had seen, convergence of the entire sequence is forced when there is a continuous dependence $\theta'' = \theta''(\theta')$. 
Even when this is too optimistic,
as the portion of missing data is likely to be smaller than the portion of observed ones, one
still gets $\theta = (\theta',\theta'',\theta''')$, where $\theta'$
is the accurate parameter, $\theta'' = \theta''(\theta')$ a portion of the spare parameter actually dependent on  $\theta'$, and therefore forced to converge,   with
$\theta'''$ 
fewer remaining spare coordinates which require additional conditions to converge.

\begin{remark}
Dimension reduction due to affine constraints on $\Theta$ (Lemma \ref{lemma1}) is not critical for the question of convergence of the EM parameter sequence $\theta^{(k)}$. 
This is different for dimension reduction due to non-minimality of the  sufficient  statistic (Proposition \ref{prop4}).
\end{remark}

\section{Alternating Bregman projections}
\label{sect_alter}

This section finds more instances where convergence of the $\theta^{(k)}$ can be guaranteed, by matching
the EM algorithm with the {\it em}-algorithm of \cite{amari}.
We consider the case of missing data from a constrained exponential family assumed minimal, where
$T(x) = (T_1(x),T_2(x)) = (y,z)$ with $y$ observed and $z$ hidden, and  we partition  $\theta = (\theta_y,\theta_z)$ 
accordingly. Minimality of
the complete data family assures that $\eta = \nabla \psi(\theta) = \mathbb E_\theta[T(x)]$ is a diffeomorphism from $G$ to $G^*$, with inverse $(\nabla \psi)^{-1} = \nabla \psi^*$,
and we may therefore
work with the expectation parameter
$\eta$. Partitioning $\eta = (\eta_y,\eta_z) = (\nabla_{\theta_y}\psi(\theta),\nabla_{\theta_z}\psi(\theta))$ in the same way, we define the data set
as 
$$
D = \{\vartheta\in \Theta: \vartheta = \nabla \psi^*(\eta), \eta_y =  y\},
$$ 
where we note  data parameters as $\vartheta\in D$, keeping $\theta\in M$ for model parameters.
We have 
\begin{lemma}
{\bf (Amari \cite{amari})}.
Let $\theta\in \Theta$. Then the right Bregman projection $\vartheta=\vec{P}_D(\theta)$ on the data set $D$ is unique, satisfies $\vartheta_z=\theta_z$,
and $\mathbb E_\theta[z | y] = \mathbb E_\vartheta[z | y]$. 
\end{lemma}

In information geometry $\vartheta_e = \vec{P}_D(\theta)$ is called the {\it e}-step from $\theta\in M$. 
It turns out that the E step from $\theta\in M$ can  also be represented as
a point $\vartheta_E \in D$ in the data set, and it generates the next M step as a left  Bregman projection
$\theta^+ \in \cev{P}_M(\vartheta_E)$.     

\begin{lemma}
{\bf (Amari \cite{amari})}. Let $\eta_E = \nabla \psi(\vartheta_E)$, $\eta_e = \nabla \psi(\vartheta_e)$ be the expectation parameters of E and {\it e}-step from
$\theta\in M$. Then $\eta_E = (y,\mathbb E_{\vartheta_e}[z|y])\in \nabla \psi(D)$ and $\eta_e = (y,\mathbb E_{\vartheta_e}[z])\in \nabla \psi(D)$.
\end{lemma}

The question when E step and {\it e}-step coincide is answered by the following:

\begin{lemma}
{\bf (Amari \cite{amari})}.  E step and {\it e}-step from $\theta$ coincide iff $\mathbb E_{\vartheta_e}[z | y] = \mathbb E_{\vartheta_e}[z]$ for $\vartheta_e = \vec{P}_D(\theta)$. If this is true all along, then
{\rm EM} and {\it em}-algorithm generate the same  iterates.
\end{lemma}

When Amari's condition  $\mathbb E(Z|Y=y) = \mathbb E(Z)$ is satisfied, we say that $Z$ is unpredictable based on knowledge of $Y$. This is a property settled between the
stronger independence (of $Z,Y$) and the weaker uncorrelatedness (cov$(Z,Y)=0$).

Let Amari's unpredictability condition be satisfied.
Then the E step is the right Bregman projection of the iterate $\theta^{(k)}\in M$ onto the data set,
$\vartheta^{(k+1)} = \vec{P}_D(\theta^{(k)})$, while the M step is the left Bregman projection
of the E step iterate $\vartheta^{(k)}$ onto the model set $M$, that is, $\theta^{(k)} \in \cev{P}_M(\vartheta^{(k)})$, 
the Bregman distance being the one induced by the log-normalizer $\psi$
of the complete data family $p(x,\theta)$. As in \cite{bregman}, we visualize this by
a building block diagram
\begin{equation}
\label{building}
\vartheta\underset{m}{\overset{l}{\longrightarrow}}\theta \underset{e}{\overset{r}{\longrightarrow}} \vartheta^{+}\underset{m}{\overset{l}{\longrightarrow}}\theta^+
\qquad\quad
\vartheta^{(k)}\underset{m}{\overset{l}{\longrightarrow}}\theta^{(k)} \underset{e}{\overset{r}{\longrightarrow}} \vartheta^{(k+1)}\underset{m}{\overset{l}{\longrightarrow}}\theta^{(k+1)}
\end{equation}

\begin{theorem}
{\bf (Convergence via information  geometry)}. Let Amari's condition be satisfied.
Then
the constrained {\rm EM} algorithm for missing data from an exponential family generates sequences $\theta^{(k)}\in M$, $\vartheta^{(k)}\in D$ of alternating Bregman projections between
$D$ and $M$, where $\theta^{(k)} \in M$ is generated by the {\rm M} step $\theta^{(k)} \in \cev{P}_M(\vartheta^{(k})$, $\vartheta^{(k)} = \vec{P}_D(\theta^{(k-1)})$ the {\rm E} step. Assume $p(x,\theta)$ is minimal, $\psi,M$ are definable,  and  suppose $M \subset G$ is bounded. Then:
\begin{enumerate}
\item[\rm (a)]
The sequence $\vartheta^{(k)}$ converges to some $\vartheta^{*}\in D$. 
\item[\rm (b)] Every accumulation point $\theta^*$ of the sequence $\theta^{(k)}$ solves the 
constrained incomplete data MLE problem, and satisfies
$\vartheta^{*} = \vec{P}_D(\theta^*)$, $\theta^* \in \cev{P}_M(\vartheta^{*})$. 
There exists $z^*$  such that $\log p(y,z^*,\vartheta^*) = \mathbb E_{\theta^*}(\log p(y,z,\theta) | y)$ is the same for every $\theta^*$.
\item[\rm (c)] Suppose $D \cap M\not=\emptyset$. Then
$D \cap M$  has a neighborhood $U$ such that any {\rm EM} sequence which 
enters $U$ converges to a point $\theta^*\in D \cap M$, i.e., $\theta^{(k)} \to \theta^*$ and $\vartheta^{(k)} \to \theta^*=\vartheta^{*}$.
\item[\rm (d)]  
The sequence $\theta'^{(k)}$ of accurate parameters converges to some $\theta'^*$, the conditional distributions converge weakly, and every $\theta^*$ gives rise to the same $P\theta^*=\theta'^*$,
$\vartheta^{*} = \vec{P}_D(\theta^*)$.
\item[\rm (e)]
When $M(\theta'^*)$ is convex, then the sequence $\theta^{(k)}$ converges to some $\theta^*\in M$ with
$\vec{P}_D(\theta^*)=\vartheta^{*}$ and $\cev{P}_M(\vartheta^{*})=\theta^*$.
\item[\rm (f)] 
Suppose the {\rm EM} instance is unconstrained. Then every sequence
$\theta^{(k)}$ with interiority converges.
\end{enumerate}
\end{theorem}

\begin{proof}
Case (a). We adopt the notation (1.4) from \cite{bregman}, where $a_k \stackrel{l}{\rightarrow} b_k \stackrel{r}{\rightarrow} a_{k+1}$ means left and right projections.  
Matching with
(\ref{building}) gives
$A=D$ and $B=M$, 
$\vartheta^{k}=a_k$ and $\theta^{(k)}=b_k$, making the results of \cite{bregman} accessible.  Then  \cite[Thm. 8.1]{bregman}
gives convergence of the $a_{k} =  \vartheta^{(k)}$. Here the $lr$-angle condition \cite[Def. 8.1]{bregman}  follows from definability
of $L,M$ and $\psi$. The $lr$-three-point inequality follows via \cite[Prop. 6.4]{bregman} from convexity of $\nabla \psi(A) = \nabla \psi(D)$,
which holds because $\nabla \psi(D)$ is the intersection of the affine
space $L=\{\eta: \eta_y =  y\}$ with im$(\nabla \psi)$.
That proves (a).

Case (b). This is a general fact using that under Amari's condition  the EM algorithm coincides with alternating Bregman projections between $D$ and $M$ as given above.

Case (c). This uses \cite[Cor. 7.4]{bregman}, which gives an even stronger statement using prox-regularity (see also Proposition \ref{prop6}).
Case (d)  is Theorem \ref{theorem2}.  
Case (e) is Theorem \ref{theorem3}. 

Case (f). The last part is when the constraint is $\xoverline{\Theta}$ and the sequence $\theta^{(k)}$ is
bounded and together with its accumulation points stays in $G$. Then we may
find a closed bounded convex definable set $M\subset G$ such that $\theta^{(k)}$
alternates between $D$ and $M$ (see \cite[Sect. 2.2]{bregman}). 
Then both projections are unique, 
$\cev{P}_M$ because $M$ is convex, and $\vec{P}_D$ because
$\nabla \psi(D) = L\cap {\rm im}(\nabla \psi)$ as the intersection of an affine subspace with
im$(\nabla \psi)$ is also convex. 
Convexity of $M$ and $\nabla\psi(D)$ also guarantees that the
$rl$- and $lr$-three-point inequalities are satisfied (see \cite[Prop. 6.4]{bregman}).
Definability of $\psi$ implies definability of $\Theta$ and
$\xoverline{\Theta}$, hence of $M$, but also definability of $\nabla \psi$, 
hence of im$(\nabla\psi)$, and since $L$ as an affine subspace is algebraic,
we get definability of $\nabla \psi(D)$.
Definability of $D$ now follows because $D$ is the image of the definable 
set $\nabla \psi(D)$ under the definable diffeomorphism $\nabla\psi^*$, see
 \cite{coste}. 
In consequence, both $rl$- and $lr$-angle conditions are satisfied
 (see \cite[Prop. 5.3]{bregman} for $rl$ and using duality \cite[Sect. 5.2]{bregman} for $lr$).

Now convergence of the sequence $b_k = \theta^{(k)}$ follows from  \cite[Thm. 7.1]{bregman}, 
while convergence of the sequence $a_k = \vartheta^{(k)}$ follows from \cite[Thm. 8.1]{bregman}.
In general the sequences converge to a gap $(\vartheta^*,\theta^*)$, that is, $\vartheta^*=\vec{P}_D(\theta^*)$, $\theta^* = \cev{P}_M(\vartheta^*)$,
possibly with $\vartheta^*\not=\theta^*$.
\hfill $\square$
\end{proof}

\begin{remark}
1)
In information geometry \cite{amari,amari_book},  $\cev{P}_M$ is called the $m$-projection,
$\vec{P}_D$  the $e$-projection. $m$-geodesics, or perpendiculars to $M$ at a left-projected point, are curved in $\theta$-coordinates, 
while $e$-geodesics, or perpendiculars to $D$ at a right-projected point,  are straight in $\theta$-coordinates, (see \cite{bregman,shawn}).
A set $M$ is $m$-flat if the left-projection
onto $M$ is unique, while a set $D$ is $e$-flat if the right-projection onto $D$ is unique. Using \cite{shawn}, and assuming $\psi$ is 1-coercive,
$m$-flat is equivalent to $M$ convex, while $e$-flat is equivalent to $\nabla \psi(D)$ convex. 

2) Alternatively, Legendreness and strict convexity of $D(x,\cdot)$ also assure
uniqueness of $\vec{P}_C$ for $C$ convex (see \cite{forward}), so here convex sets are also
$e$-flat. However, this is less useful in the present setting, where it is $\nabla \psi(D)$ which is convex, not $D$.

3) In the information geometry literature statement (f) has been made repeatedly, but without the hypothesis of definability of $\psi$. We are aware of a couple of published 
incorrect proofs.
Our own proof requires the K\L-condition,  and one would of course like to know whether this can be avoided. Note that
we can treat the case $D \cap M=\emptyset$.

4) Case (f) can also be derived from Theorem \ref{theorem3}.

5) Case (f) may have local minima or convergence to saddle points, which confirms that even this simplest instance of EM is not 
from the realm of convexity despite $-\log q(y,\cdot)$ and $-\log p(x^+,\cdot)$ being convex. 
A simple case with two points in $D\cap M$ and a non-zero gap between $D$ and $M$
is given in Example \ref{ex_arc}.

6) The case $\phi(s) = s^{1/2}$ in the K\L-condition gives linear convergence. This occurs for instance in the neighborhood
of a point $\theta^\sharp\in D \cap M$ where $D,M$ intersect transversally. (For the definition of transversal intersection in the Bregman sense
see \cite{bregman}).
\end{remark}

\begin{remark}
Instead of {\it e}- and {\it m}-flatness it is
preferable
to requests uniqueness of the projections only for points close enough.  That leads to the concept of prox-regularity, or positive reach,
which is used in \cite{bregman}. The significant difference is that positive reach is invariant under duality, i.e.,
$D$ has positive reach iff $\nabla \psi(D)$ has, and the same for $M$, while auto-duality fails for $e$- and $m$-flatness. 
A sample result is:
\end{remark}

\begin{proposition}
\label{prop6}
Consider the missing data case
under the hypotheses of Theorem {\rm \ref{theorem2}}. Suppose the set $M$ has positive left Bregman reach $r^* > 0$ at the accumulation points $\theta^*$
of the sequence $\theta^{(k)}$. Let $\vartheta^*=\vec{P}_D(\theta^*)$ and
suppose $\min_{\theta \in M}D(\theta,\vartheta^{*}) <\frac{1}{2} r^{*2}$. Then the sequence $\theta^{(k)}$ converges. 
\end{proposition}

\begin{proof}
We know that the data set sequence $\vartheta^{(k)}$ converges to $\vartheta^{*}$ and $\theta^{(k)}\in \cev{P}_M(\vartheta^{(k)})$. 
But $\cev{P}_M(\vartheta^{*}) = \theta^{\diamond}$ is singleton due to the hypothesis on left Bregman reach of $M$. That clearly
implies $\theta^{(k+1)} \to \theta^\diamond$. 
\hfill $\square$
\end{proof}

\section{More general families}
\label{sect_beyond}
Several of the arguments used for exponential families can be extended to a more general setting.
We consider the case where $d\mathbb P_\theta = p(T(\cdot),\theta) d\mu$, and for the conditional family,
$d\mathbb P^{x|y}_\theta= k(T(\cdot)|y,\theta) d\mu_y$ for a sufficient statistic $T(x)$. Suppose there is an orthogonal
change of coordinates $Q\theta = (\theta',\theta'')$ 
and a possibly non-linear reduction to a minimal sufficient statistic $T'(x)$ such that the conditional family has the equivalent representation
$$
d\mathbb P^{x|y}_\theta =d\mathbb P^{x|y}_{\theta'} = k'(T'(x)|y,\theta') d\mu_y(x),
$$
depending only on $\theta'$. Let us consider the following
property extending affine independence of the $T_j$ in the case of exponential families (Definition \ref{def2}):
\begin{eqnarray}
\label{minimal}
\begin{array}{l}
\mbox{In an affine-minimal representation ${k'}(T'(x)|y,{\theta'})$ of $\mathbb P^{x|y}_\theta$ there exists}\\ 
\mbox{no $v\not=0$ with $v \cdot \nabla_{{\theta'}} \log {k'}({T'}(x)|y,{\theta'})=0$ for $\mu_y$-almost all $x\in h^{-1}(y)$.}
\end{array}
\end{eqnarray} 

\begin{proposition}
Suppose $\mathbb P^{x|y}$ has the above property, and  let $k'(T'(x)|y,\theta')d\mu_y(x)$ be an  affine-minimal representation
of $\mathbb P^{x|y}_\theta$.
Then $\nabla^2_{22} \bar{K}_y({\theta'}||{\theta'}) \succ 0$.
\end{proposition}

\begin{proof}
Let $k(\cdot|y,\theta) d\mu_y$ be affine-minimal for the ease of notation.
Since $
\nabla^2_{22} K_y(\theta||\theta)= 
\mathbb E_\theta[ \nabla_\theta k(\cdot|y,\theta) \nabla_\theta k(\cdot|y,\theta)^T | y] \succeq 0$
by (\ref{stimmts}),
$v \cdot \nabla^2_{22} K_y(\theta||\theta) v=0$
implies 
$\mathbb E_\theta \left[ |v \cdot \nabla_\theta k(\cdot|y,\theta)|^2|y\right]=0$,
hence
$v\cdot \nabla_\theta \log k(x|y,\theta)=0$ $\mu_y$-a.s., which by minimality implies $v=0$.
Therefore $\nabla_{22}^2 K_y(\theta||\theta) \succ 0$.
\hfill $\square$
\end{proof}

This means we can get a situation as previously found for the exponential family. There exists an orthogonal $n\times n$-matrix $Q$
such that
$$
Q\nabla^2_{22}K_y(\theta||\theta)Q^T =   
 \left[\begin{array}{c|c} \nabla^2_{22} \bar{K}_y(\theta'||\theta') & 0 \\ \hline
0 & 0 \end{array}\right], \quad Q \theta = \left[ \begin{array}{c} {\theta'}\\\ \!\!{\theta''}\end{array}\right], \quad  
\nabla^2_{22}\bar{K}_y(\theta' ||\theta')  \succ 0,
$$
where $\theta'$ are the accurate parameters remaining after removing the affine dependence in $\theta$. 
In consequence, we can again prove convergence of the accurate parameter sequence $\theta'^{(k)}$ using the partial convergence theorem.

\begin{theorem}
\label{thm6}

Suppose $-\log q(y,\cdot)+ i_M$ satisfies the K\L-inequality and  $M \subset G$ is bounded.
Suppose the conditional family has the minimality property above. Let $\theta^{(k)}$ be generated by the constrained EM algorithm. Then
the sequence $\theta'^{(k)}$ of accurate parameters converges.
\end{theorem}

The last step is to consider the split program $(P_\infty)$. We have the following
result, which uses the proof of Proposition \ref{principle}.

\begin{proposition}
Under the hypotheses of Theorem {\rm \ref{thm6}}.
Suppose the conditional expected Fisher information matrix 
$\mathcal I_c(\theta,y)$ of complete data given $y$ is positive definite
on $G$. Then we have $-\nabla^2_{\theta'' \theta''} \log q(y,\theta',\theta'') \succ 0$ for fixed $\theta'$.
\end{proposition}

This allows again to upgrade convergence of $\theta'^{(k)}$ to convergence of the full parameter sequence, e.g. for
$M(\theta	')$ convex,  when the 
complete data family has a full rank expected conditional Fisher information of complete data, given $y$, i.e.,
$\mathcal I_c(\theta,y) \succ 0$.

\section{Examples}
\label{sect_examples}
In this section we discuss several limiting examples.

\begin{example}
\label{ex_missing}
We consider the missing data  case $x=(y,z)$, where, $h:(y,z) \mapsto y$ with $y$ observed and $z$ hidden,  and with the special statistic
$T(x) = x$. Let $d\mu = m(y,z) d\mu_Y \otimes d\mu_Z$, where $\mu_Y \otimes \mu_Z$ is a product measure on $X = Y \times Z$.
We have $p(x,\theta) = e^{ \theta_y\cdot y +  \theta_z\cdot z - \psi(\theta_y,\theta_z)}$.
In the notation of the disintegration, $\nu = \mu_Y$ and $d\mu_y = m(y,z) d(\delta_{\{y\}} \otimes \mu_Z)$  for  $y\in Y$, because
\begin{align*}
\int_{Y\times Z} f(y,z) d\mu(y,z) &= \int_Y \left[ \int_Z f(y,z) m(y,z) d\mu_Z(z)\right] d\mu_Y(y)\\
&= \int_Y \left[  \int_{\{y\} \times Z} f(y,z) m(y,z) d(\delta_{\{y\}} \otimes \mu_Z)(y,z) \right] d\mu_Y(y)\\
&= \int_Y \left[\int_{h^{-1}(y)} f(y,z) d\mu_y(y,z) \right] d\nu(y).
\end{align*}
That gives
\begin{align*}
\psi_y(\theta) &= \log \int_{h^{-1}(y)} e^{\theta_y\cdot y} e^{\theta_z\cdot z} d\mu_y(y,z) \\
&= \log \int_{\{y\}\times Z} e^{\theta_y\cdot y} e^{\theta_z\cdot z} m(y,z) d(\delta_{\{y\}} \otimes \mu_Z)(y,z)\\
&={\theta_y\cdot y} + \log \int_Z e^{\theta_z\cdot z} m(y,z) d\mu_Z(z).
\end{align*}
Here $T(h^{-1}(y)) = \{y\} \times Z$, so that we need  dimension reduction in Proposition \ref{prop4}.
We obtain
\begin{align}
\label{wenig}
\begin{split}
\log k(y,z|y,\theta_y,\theta_z)
&= {\theta_y\cdot y} + {\theta_z\cdot z} - {\theta_y\cdot y} -\log \int_Z e^{\theta_z\cdot z} m(y,z) d\mu_Z(z)\\
&=: \theta_z\cdot z - \chi_y(\theta_z) =: \log \bar{k}(z|y,\theta_z) 
\end{split}
\end{align}
on defining $\chi_y(\theta_z) = \log \int_Z e^{\theta_z\cdot z} m(y,z) d\mu_Z(z)$.
Therefore the conditional family is
$$
k(y,z|y,\theta_y,\theta_z) d\mu_y(y,z) = \bar{k}(z|y,\theta_z) d\mu_Z(z)
$$
and the accurate parameter is $\theta_z$, at least when $z$ is minimal for the conditional family, the spare parameter being $\theta_y$. 
That means even under definability 
we can only expect convergence of the  $\theta_z$-part of $\theta$. 
For
convergence of the $\theta_y$-part we must rely
on the split program $(P_\infty)$.

\begin{example}
{\bf (Continued...)}. We compare the domains of objective and regularizer in (\ref{missing}).
From $q(y,\theta_y,\theta_z) = \int_Z e^{\theta_y\cdot y} e^{\theta_z\cdot z}
e^{-\psi(\theta_y,\theta_z)} m(y,z) d\mu_Z(z)
= e^{\theta_y\cdot y-\psi(\theta_y,\theta_z)}\int_Z e^{\theta_z\cdot z} m(y,z) d\mu_Z(z)$
follows dom$\,q(y,\cdot) = {\rm dom}\, p(x,\cdot) = {\rm dom}(\psi)=
\Theta =   \{(\theta_y,\theta_z): \int_Z e^{\theta_y\cdot y} e^{\theta_z\cdot z} m(y,z) d\mu_y \otimes \mu_Z < \infty\}$, while dom$\,\psi_y = Y \times \{\theta_z: \int_Z e^{\theta_z\cdot z} m(y,z) d\mu_Z(z)<\infty\}= {\rm dom}\, \bar{K}_y(\theta^{(k)}|| \cdot)$ is larger. Nonetheless, $\Theta$ and {\rm dom}$\,\psi_y$
have common boundary points.
\end{example}

\end{example}

\begin{example}
\label{ex_amari}
Taken from \cite{amari}. Consider an independent normal sample
$x_1,\dots,x_N$ with statistic
$T(x) = (T_1(x),T_2(x))=(\sum_{i=1}^N x_i/N, \sum_{i=1}^N x_i^2/N)$,
and suppose  $y=T_1(x)$ is observed, while $z=T_2(x)$ is hidden. 
 The corresponding exponential family 
$p(x,\mu,\sigma^2)$ is written
as $p(x,\theta)$ with $\theta = (\theta_1,\theta_2) = (N\mu/\sigma^2,-N/2\sigma^2)$, hence $\mu=-\theta_1/2\theta_2$, $\sigma^2=-N/2\theta_2$,
$\Theta = G=\mathbb R \times (-\infty,0)$.

Specializing to $N=2$ for simplicity, we have $\psi(\theta_1,\theta_2) = -\theta_1^2/4\theta_2-\log(-\theta_2)$ defined on $G$,
with Legendre transform
$\psi^*(\eta_1,\eta_2) = -1-\log(\eta_2-\eta_1^2)$ defined on $G^*=\{\eta:\eta_2 > \eta_1^2\}$,  gradient, inverse gradient and expectation parameter being
$$
\nabla \psi(\theta) = \begin{bmatrix} -\frac{\theta_1}{2\theta_2} \\ \frac{\theta_1^2}{4\theta_2^2} - \frac{1}{\theta_2} \end{bmatrix}
=: \begin{bmatrix} \eta_1\\\eta_2\end{bmatrix}
\qquad
\theta_1 = -\frac{2\eta_1}{\eta_1^2-\eta_2}, \theta_2 = \frac{1}{\eta_1^2-\eta_2}, \qquad
\nabla \psi^*(\eta) = \begin{bmatrix}
-\frac{2\eta_1}{\eta_1^2-\eta_2} \\ \frac{1}{\eta_1^2-\eta_2} 
\end{bmatrix}
$$
and in the original coordinates
$\eta_1 = \mu$, $\eta_2 = \mu^2+\sigma^2$.
We have $h(x_1,x_2) = y = \frac{x_1+x_2}{2}$. The family is two-dimensional,
but $T(h^{-1}(y)) = (y, \frac{1}{2}(x_1^2 + (2y-x_1)^2))$ is included in the one-dimensional affine
space $T_1=y$ of $(T_1,T_2)$. Hence we need parameter reduction. The conditional family is
\begin{align*}
k(x_1,x_2|y,\theta) &= C(\sigma) \exp \{- [(x_1-\mu)^2+(x_2-\mu)^2]/{2\sigma^2}\} \big/ \exp\{ -(y-\mu)^2/\sigma^2\} \\
&=  \exp\{-{(x_1-y)^2}/{\sigma^2}\}/2\sqrt{\pi} \sigma =: \bar{k}(x_1|y,\theta_2) \\
&= \exp\left\{ - (x_1^2-2yx_1)/\sigma^2 - \left[   {y^2}/{\sigma^2}- \textstyle\frac{1}{2} \log {\sigma^{-2}}\right]\right\} /{2\sqrt{\pi}}\\
&= \exp\left\{ \theta_2 \cdot (x_1^2-2yx_1) - \left[ -y^2\theta_2 - \textstyle\frac{1}{2} \log(-\theta_2)  \right]\right\} /{2\sqrt{\pi}}
\end{align*}
the reduced family depending only on $\theta'=\theta_2$, which is the accurate parameter.  
The incomplete data family is obtained as follows: Since $\mathbb E_\theta [(x_1+x_2)/2] = \mu$ and $\mathbb V_\theta[(x_1+x_2)/2] = \sigma^2/2$,
we have $q(y,\theta) \sim N(\mu,\sigma^2/2)$.

Now we consider the M step. We discuss two cases, constrained and unconstrained.
In the first scenario
a constraint is introduced in the form of a curved exponential family, namely $\mu^2 = \sigma^2$, so that the model family is $N(\mu,\mu^2)$. 
In natural parameters this is
$M = \{\theta: \theta_1^2 = -4\theta_2\}$. Here Theorem \ref{theorem2} assures convergence of the accurate parameter
sequence $\theta_2^{(k)}$. But the constraint gives
$\theta_1^2$ as a function of $\theta_2$, so $\theta_1^{(k)}$ converges, too.

In the unconstrained case the parameter $\theta_2^{(k)}$ converges, hence so does $\sigma^{(k)2}$,
while the parameter $\theta_1$ is free. We have to consider 
the M step program $(P_k)$, which is
$
\min_{\theta_1} -\log q(y,\theta_1,\theta_2^{(k)}). 
$
Since $-\log q(y,\theta_1,\theta_2^{(k)}) = (y-\mu)^2/\sigma^{(k)2} + \log C(\sigma^{(k)})$, the solution is
always $\mu=y= \mu^{(k)}$, which implies $\theta_1^{(k)} = 2y/\sigma^{(k)2} = -2y\theta_2^{(k)}$, which again converges. This is interesting, as in
the first place the available information does not seem sufficient to estimate $\mu$ and
$\sigma^2$.
\end{example}

\begin{example}{\bf (Continued...)}.
As observed in \cite{amari}, in this example a difference between the EM algorithm and the {\it em}-algorithm occurs. 
In our present terminology this is due to the fact that
$\mathbb E_\vartheta[(x_1^2+x_2^2)/2] = \mu^2+\sigma^2=y^2+\sigma^2$ and $\mathbb E_\vartheta[(x_1^2+x_2^2)/2 | y] = y^2 + \sigma^2 /2$
are different. 

This discrepancy can be easily remedied by
letting $T_2(x) = \sum_{i=1}^N (x_i-\bar{x})^2/(N-1)$, as then $T_1,T_2$ are independent, so that Amari's condition is satisfied.
\end{example}

\begin{example}
\label{ex_arc}
Consider the Kullback-Leibler distance in $\mathbb R^2_+$ given as
$K(a||b) = \sum_{i=1}^2 a_i \log \frac{a_i}{b_i} - a_i + b_i$, which is a special Bregman distance induced by $\psi(x) = \sum_{i=1}^2 x_i \log x_i - x_i$. 
Choose two points $p,q\in (-\infty,0)^2$ and compute
$a = \exp(p)$, $b=\exp(q)$. Let $p_t = tp + (1-t)q$, $t\in [0,1]$ be the points on the segment $[p,q]$, then
$a_t = \exp(p_t) = \exp(tp)\exp((1-t)q)$ forms a curve in $(0,1)^2$, which in general is not straight. Let $A = \{a_t: t \in [0,1]\}$. 
Then $\nabla \psi(A) = \log(A) = [p,q]$ the segment, hence $A$ is {\it e}-flat. In other words, Kullback-Leibler right projections
on $A$ are unique. Now let $B= [a,b]$ be the segment joining $a,b$. Then $B$ is convex, hence Kullback-Leibler left projections on $B$
are unique, and $B$ is {\it m}-flat. We have $A \cap B = \{a,b\}$ in the case where $A$ is curved. This means Bregman projections
$\vec{P}_A$, $\cev{P}_B$ are unique, and $\vec{P}_A \circ \cev{P}_B$  has two points of attraction.
However, starting with a point $a_t$ somewhere in between, one can see iterates either go left, or go right
toward $a$ or $b$. Except a point pair $a',b'$ which satisfies $a' = \vec{P}_A(b')$,
$b'=\cev{P}_B(a')$, building a gap or fixed point pair. If started at $a'$, the alternating projection method will remain at the fixed point pair. 
\end{example}

\begin{example}
{(\bf  Failure of convergence for curved family).}
    We consider independent gaussian random variables $X_1,X_2,X_3$
    with $\mathbb E(X_i) = \mu_i$ and $\mathbb V(X_i) = 1$. 
    We want to estimate $\mu=(\mu_1,\mu_2,\mu_3)$ based on an 
    observation $y$ of $Y=X_3$, where it is  
   assumed that $\mu_3 = f(\mu_1,\mu_2)$ with a known $f$. Concerning the statistic
   that means  $x=(x_1,x_2,x_3) = (z_1,z_2,y)$ with $y$ observed and $z$ hidden.     
   This
    gives the model set $M=\{\mu: \mu_3 =f(\mu_1,\mu_2)\}$, hence the family is curved.
   Suppose $y=0$.
   Then the E step  $x^+=\mathbb E(X|Y=y,\mu)$ yields $x^+_1=\mu_1,x^+_2=\mu_2$, $x^+_3=0$.
The M step is $\mu^+ \in \argmin_{\mu\in M} -\log p(x^+,\mu)$, where
$$
p(x^+,\mu) = C \exp\{ -\textstyle\frac{1}{2}(x_1^+-\mu_1)^2-\frac{1}{2}(x_2^+-\mu_2)^2-\frac{1}{2}(x_3^+-\mu_3)^2\}.
$$
Hence the E step
is  orthogonal projection of $\mu\in \mathbb R^3$ onto the data set $D= \{x_3=0\}$. 
   The M step is orthogonal projection of $x^+=(x_1^+,x_2^+,0)$ onto the model set $M={\rm graph}(f)$. 
The result of this projection is $\mu^+\in M$, and then the procedure is repeated.
We therefore have a case, where EM and {\it em}-algorithm coincide. 

Altogether the method is now the alternating projection method between the sets
$A = \{(x_1,x_2,x_3): x_3=0\}$ and $B = \{(x_1,x_2,x_3): x_3 = f(x_1,x_2)\}$, and this can readily be extended to $x\in \mathbb R^n$,
where
$A=\{(x,0): x\in \mathbb R^n\}$ and $B= \{(x,f(x)): x\in \mathbb R^n\}$ the graph of a function
$f:\mathbb R^n\to \mathbb R$.  Assume $f(x) \geq 0$. 
Then 
\begin{equation}
\label{pr_pp}
P_A(x_k,f(x_k)) = (x_k,0), \quad
(x_k,f(x_k)) \in P_B(x_{k-1},0) \; \mbox{ iff }\;  x_{k-1} = x_k + f(x_k) \nabla f(x_k).
\end{equation}
Infinitesimally, this method follows steepest ascent
backwards. 

Going back to $n=2$, 
we let $B=M$  be the graph of the mexican hat function \cite{absil} on $x_1^2+x_2^2 \leq 1$.
Similar to the argument given for steepest descent with infinitesimal steps in \cite{absil},  AP with infinitesimal steps will also follow the valley of the hat downward, 
endlessly circling around and approaching the boundary curve $x_1^2+x_2^2=1$, where 
$f=0$.  
What is amiss for convergence is the  K\L-property of  $M$, which is not definable. For a picture see \cite{absil}.

    \end{example}

\begin{example}
{\bf (PPM and EM).}
Take again the situation $A=\{(x,0): x\in \mathbb R^n\}$,
$B = \{(x,f(x)): x\in \mathbb R^n\}$, where $f \geq 0$. Consider the proximal point step
\begin{equation}
\label{pr}
x_k \in \argmin \textstyle\frac{1}{2}f(x)^2 + \frac{1}{2}\|x-x_{k-1}\|^2.
\end{equation}
The necessary optimality condition is $0 = f(x_k)\nabla f(x_k) + x_k-x_{k-1}$, which is the alternating projection step (\ref{pr_pp}) above. This means,
(\ref{pr}) must be the scheme (\ref{general}), respectively, its realization in Proposition \ref{prop1}, for our example above.
Choosing $f$ such that AP fails to converge, we also produce an example, where the non-convex proximal point method (with fixed $\lambda_k=1$) fails to converge, now with objective $\frac{1}{2}f(x)^2$. This construction makes every instance of PPM with an objective bounded below a special instance of EM.
\end{example}

\end{document}